\documentclass[12pt]{article}

\usepackage{amssymb,a4}
\usepackage{amsmath,amsfonts,amssymb,amsthm, mathrsfs}
\usepackage{enumerate}

\usepackage{upgreek}
\usepackage{wasysym}
\usepackage[all,cmtip,dvipdfmx]{xy}
\usepackage{xcolor}
\definecolor{rouge}{rgb}{0.85,0.1,.4}
\definecolor{blue}{rgb}{0.1,0.2,0.9}
\definecolor{violet}{rgb}{0.7,0,0.8}

\usepackage{color}
\usepackage{colordvi}

\usepackage[
colorlinks=true,linkcolor=blue,urlcolor=violet,citecolor=rouge]{hyperref}

\newcommand{\on}{\operatorname}

\newcommand{\al}{\alpha}

\newcommand{\End}{\operatorname{End}}

\def\be{\beta}

\def\wW{\widehat{W}}

\newcommand{\la}{\lambda}
\def\leq{\leqslant}
\def\geq{\geqslant}

\def\C{{\mathbb C}}
\def\Q{{\mathbb Q}}

\def\Z{{\mathbb Z}}

\def\N{{\mathbb N}}
\def\1{{\bf 1}}

\def \End{{\rm End}}

\def \mod{{\rm mod}}

\def \<{\langle}
\def \>{\rangle}

\def \g{\frak g}
\def \fg{\widehat{\frak{g}}}
\def\fh{\widehat{\frak{h}}}
\def \n{\frak{n}}

\def\wr{\widehat{\rho}}

\def \wD{\widehat{\Delta}}

\def \wW{\widehat{W}}
\def\wD{\widehat{\Delta}}

\def \h{\mathfrak{h}}

\def\wQ{\widehat{Q}}
\numberwithin{equation}{section}

\newtheorem{theorem}{Theorem}[section]
\newtheorem{prop}[theorem]{Proposition}
\newtheorem{lem}[theorem]{Lemma}
\newtheorem{cor}[theorem]{Corollary}
\newtheorem{remark}[theorem]{Remark}
\newtheorem{conjecture}[theorem]{Conjecture}

\theoremstyle{definition}
\newtheorem{defn}[theorem]{Definition}

\begin{document}
 
  \begin{center}{Associated varieties of simple affine VOAs $L_k(sl_3)$ and $W$-algebras $W_k(sl_3,f)$}
\end{center}
\begin{center}
	{ Cuipo Jiang
		and Jingtian Song}\\
	School of Mathematical Sciences, Shanghai Jiao Tong University, Shanghai 200240, China
\end{center}

 
\begin{abstract} 
In this paper we first prove that the maximal ideal of the universal affine vertex operator algebra $V^k(sl_n)$  for $k=-n+\frac{n-1}{q}$  is generated by two singular vectors of conformal weight $3q$ if $n=3$, and by one singular vector of conformal weight $2q$ if $n\geq 4$.  We next determine the associated varieties of the simple vertex operator algebras $L_k(sl_3)$ for all the non-admissible levels $k=-3+\frac{2}{2m+1}$, $m\geq 0$.  
The varieties of the associated simple affine $W$-algebras $W_k(sl_3,f)$,  for  nilpotent elements $f$ of $sl_3$,  are also determined. 
 \end{abstract}

\section{Introduction}
Let $V$ be a vertex algebra. Recall from \cite{Zhu96} that the 
 {\em Zhu's $C_2$-algebra}  of $V$ is by definition the quotient space
$R_V=V/C_2(V)$ of $V$, 
where $C_2(V)=\on{span}_{\C}\{a_{-2}b\mid a,b\in V\}$. 
Then $R_V$ has the Poisson algebra structure given by
$$\bar a\cdot\bar b=\overline{a_{-1}b},\qquad \{\bar a,\bar b\}=\overline{a_{0}b},$$
for $a,b \in V$, where   $\bar a := a+C_2(V)$ is the image of $a$ in $R_V$.
The associated variety $X_V$ of $V$ is the reduced scheme
$X_V=\on{Specm}(R_V)$
of $R_V$ \cite{Ar12a}.
This fundamental   invariant of $V$ has been extensively studied in \cite{BFM,Zhu96,ABD04,M04,Ar12a,Ar15a,Ar15b,AM18a,ArK18, AM17,AJM21, ADFLM24,AFK24}.

Let $\g$ be a complex simple  finite-dimensional  Lie algebra and $V^k(\g)$ 
the universal affine vertex algebra  at level $k \in\C$. By PBW Theorem, $R_{V^k(\g)}$ is isomorphic to the polynomial algebra $\C[\g^*]$ over $\g^*$. 
Thus the variety $X_V$  is  the affine space $\g^*$ with Kirillov-Kostant Poisson structure.
For a graded quotient $V$
of $V^k(\g)$, the variety $X_V$ is a Poisson subscheme of $\g^*$
which is $G$-invariant and conic, where $G$ is the adjoint group of $\g$.

Denote by $L_k(\g)$ the simple quotient of $V^k(\g)$.  
Let $h$ be the Coxeter number of $\g$, $h^\vee$ the dual Coxeter number, and $r^\vee$ the lacing number of $\g$.
It has been proved in \cite{AJM21} that 
$L_k(\g)=V^k(\g)$ if and only if $X_{L_k(\g)} = \g^*$.   In \cite{GK07}, Gorelik and Kac proved 
that $V^k(\g)$
is not simple if and only if
\begin{align}
r^{\vee}(k+h^{\vee}) \in \Q_{\geq 0}\backslash \left\{\frac{1}{m}\mid m\in \Z_{\geq 1}\right\}.
\label{eGK}
\end{align}
Therefore,
\begin{align}
X_{L_k(\g)}\subsetneq \g^*\iff \text{(\ref{eGK}) holds.}
\end{align}
It is known \cite{Zhu96,DM06, Ar12a} that $X_{L_k(\g)}=\{0\}$
 if and only if  $k$ is a non-negative integer. 
For  an {\em admissible}   $L_k(\g)$ \cite{KW89},
or equivalently, for an {\em admissible} number $k$ which by definition satisfies 
\begin{align*}
k+h^{\vee}=\frac{p}{q},\quad p,q\in \Z_{\geq 1}, \ (p,q)=1,\
p\geq \begin{cases}h^{\vee}&\text{if }(r^{\vee},q)=1,\\
h&\text{if }(r^{\vee},q)\ne 1,
\end{cases}
\end{align*}
the variety $X_{L_k(\g)}$ is the closure of some nilpotent orbit in $\g$  \cite{Ar15a}.  In particular,  if $k$ is a non-degenerate admissible number, that is,
 \begin{align*}
 q\geq \begin{cases}h &\text{if }(r^{\vee},q)=1,\\
 r^{\vee}h^{\vee} &\text{if }(r^{\vee},q)=r^{\vee},
 \end{cases}
 \end{align*}
 then  $X_{L_k(\g)}$ is the nilpotent cone ${\mathcal N}(\g)$ \cite{Ar15a}.
 If $k=-h^{\vee}$, 
 it follows from \cite{FF92, FG04} that 
 $X_{L_k(\g)}$ is also the nilpotent cone.
 It was observed in  \cite{AM18a,AM18b,AM17,ADFLM24}
 that
there are   
 cases when $L_k(\g)$ is
non-admissible 
and
$X_{L_k(\g)}$ is the closure of some nilpotent orbit,
 which provide examples for the conjecture in physics \cite{XY} that,
 in view of the 4D/2D duality,
 there should be
 a large list of non-admissible simple affine vertex algebras
 whose associated varieties are the closures of some nilpotent orbits.  
 There are also cases  where $X_{L_k(\g)}$
 is neither $\g^*$ nor contained in the nilpotent cone ${\mathcal N}(\g)$.
 It was  proved in \cite{AM17}  that $X_{L_{-1}(sl_n)}$, $n\geq 4$, $X_{L_{-m}(sl_{2m})}$, $m\geq 2$, $X_{L_{2-r}(so_{2r})}$, $r\in 2\Z_++1$,  are closures of  sheets. In \cite{AFK24}, it was shown that the associated variety of a simple affine vertex algebras is contained in the closure of a Dixmier sheet when a chiralization of generalized Grothendieck's simultaneous resolution exists. 
 The variety  $X_{L_{-\frac{5}{2}}(sl_4)}$ of $L_{-\frac{5}{2}}(sl_4)$ was given in \cite{Fa22}, which turns out to be the closure of a Jordan class of $sl_4$. 
 
 \smallskip
In general,
the problem of determining $X_{L_k(\g)}$ is widely open. Up to now, the  variety $X_{L_k(\g)}$ has been completely determined for all $k\in\C$ only when $\g=sl_2$. 
In this paper, we will determine $X_{L_k(sl_3)}$ for all  non-admissible  numbers $k=-3+\frac{2}{2m+1}$, $m\geq 0$. Then together with the  results $X_{L_k(\g)}$
for the critical case and admissible cases given in \cite{FF92, FG04,Ar12a}, we will obtain a complete characterization of $X_{L_k(sl_3)}$ for all $k\in\C$.  

\smallskip
Another question we are concerned about in this paper is the characterization of the maximal ideal of $V^k(sl_n)$, where $k=-n+\frac{n-1}{q}$, $n\geq 3$, $q\geq 2$, $(n-1, q)=1$. 
 If $q=1$, that is,  $k=-1$, the maximal ideal of $V^{-1}(sl_3)$ is generated by two singular vectors of conformal weight 3 \cite{AP08, AKMPP20}, and the maximal ideal of $V^{-1}(sl_n)$ for $n\geq 4$ is generated by one singular vector of conformal weight 2 \cite{AP08, AM17}.
By  using the character formulas given in \cite{KT00}, the fact that the coefficients as values of Kazhan-Lusztig (or inverse KL) polynomials  at 1 in the formulas depend only on the corresponding Coxeter group,  the category equivalences  given in \cite{Fi06},  and the  results for the $k=-1$ case in \cite{KW18, BKK24, AP08, AM17, AKMPP20}, we prove our first main result  as follows (also see Theorem \ref{tmain3.1}).
\begin{theorem}\label{main1}  
	Let $k=-n+\frac{n-1}{q}$ with $n\geq 3$, $q\geq 2$, $(n-1,q)=1$. Then 	
$$
	\widehat{R}{\rm ch}L_k(sl_n)=
	\sum\limits_{w\in W}\sum\limits_{\tiny{\begin{gathered}\gamma\in Q \\(\gamma|\bar{\Lambda}_{n-1})\geq 0, (\gamma|\bar{\Lambda}_1))\geq 0\end{gathered}}}(-1)^{l(w)}e^{wt_{q\gamma}(k\Lambda_0+\wr)-\wr}.	
$$
	Furthermore, we have 
	\begin{enumerate}
		\item If $k=-3+\frac{2}{2m+1}$ with $m\geq 1$, then  the maximal ideal of $V^k(sl_3)$ is generated by  two singular vectors $v^1$ and $v^2$ of weights $k\Lambda_0-3(2m+1)\delta+2\alpha_1+\alpha_2$ and $k\Lambda_0-3(2m+1)\delta+\alpha_1+2\alpha_2$, respectively. 
		
		\item If $k=-n+\frac{n-1}{q}$ with $n\geq 4$, $q\geq 2$, and $(q,n-1)=1$, then the maximal ideal of $V^k(sl_n)$ is generated by one singular vector $w^1$ of weight $k\Lambda_0-2q\delta+\al_1+2\al_2+\cdots+2\al_{n-2}+\al_{n-1}$. 
	\end{enumerate}
\end{theorem}

\smallskip
Let $\theta$ be the highest root of $\g$. 
Our second main result is (also see Theorem \ref{variety1}):
\begin{theorem}\label{main2}
	The variety $X_{L_{-1}(sl_3)}$ is the closure of the Dixmier sheet   ${\mathbb S}_{\frak{l}_1}=G.{\mathbb C}^*\lambda$, i.e.,
	$$
	X_{L_{-1}(sl_3)}=\overline{{\mathbb S}_{\frak l_1}}= \overline{G.{\mathbb C}^*\lambda}=G.{\mathbb C}^*\lambda\cup\overline{{\mathcal O}_{\min}},
	$$
	where ${\frak l}_1=\frak{h}\oplus \C e_{\theta}\oplus \C f_{\theta}$, $\lambda=h_1-h_2$, and ${\mathcal O}_{\min}$ is the minimal non-zero nilpotent orbit of $sl_3$.  In particular, $X_{L_{-1}(sl_3)}$ is irreducible and 
	$
	\dim X_{L_{-1}(sl_3)}=5
	$.
\end{theorem}
\begin{remark}
The result of Theorem \ref{main2} is also covered in \cite{AFK24} in a different and more conceptual way.
\end{remark}
For $k=-3+\frac{2}{2m+1}$, $m\geq 1$, we have our third main result as follows (see also Theorem \ref{t3.2.1}).
\begin{theorem}\label{main3}
For $k=-3+\frac{2}{2m+1}$, $m\geq 1$, the variety
	$$X_{L_k(sl_3)}= \overline{G.({\mathbb C}^*(h_1-h_2)+f_{\theta})}.$$
\end{theorem}
The proof of  Theorem \ref{main3} is much more complicated than the $k=-1$ case. The main difficulty is that although we could prove that the maximal ideal of  $L_k(sl_3)$ is generated by two singular vectors $v^1$ and $v^2$ of conformal weight $3(2m+1)$ by Theorem \ref{main1}, it is almost impossible to write these  singular vectors out  as sums of PBW terms for general $k=-3+\frac{2}{2m+1}$, $m\geq 1$.  This means that we could not use their structure, as what we do with the  $k=-1$ case, to obtain the zero locus of the  graded ideal $\overline{{\mathcal I}_k}$ of the $C_2$-algebra $R_{V^k(sl_3)}$ determined by these two singular vectors. Our main idea is to analyze coefficients of some critical leading terms.  Then  further prove that 
$X_{L_k(sl_3)}= \overline{G.({\mathbb C}^*(h_1-h_2)+f_{\theta})}$.  Our method here provides a possible way to determine varieties of a simple affine vertex operator algebra without knowing the exact structure of its singular vectors.  We will study the associated varieties of the simple vertex operator algebras $V^k(sl_n)$ for $n\geq 4$ and general non-admissible numbers $k$ in subsequent papers.

\smallskip
Representations of $L_k(sl_3)$ have been extensively studied in \cite{AP08, Ad16, AFR17,  Cr18, FLLZ18, AdMi21, AdK21, AKR21,  FKR21, FR22,  KRW22, AMP23, CRR21, ACG21, AKR24}, etc. 
We  have the following conjecture for $L_k(sl_3)$ (also see Conjecture \ref{conj1}).
\begin{conjecture} Let $k=-3+\frac{2}{2m+1}$ with $m\geq 0$. Then
$$
\{L(t\bar{\Lambda}_1-\frac{2i}{3}\bar{\Lambda}_2), L(t\bar{\Lambda}_2-\frac{2i}{3}\bar{\Lambda}_1),  L(\bar{\Lambda}_1-(t+\frac{2i}{3}+1)\bar{\Lambda}_2),  ~t\in\C, i=0,1,\cdots, 2m\}$$
provides a complete list of irreducible $A(L_k(sl_3))$-modules in the category ${\mathcal O}$.
\end{conjecture}

When $m=0$, the conjecture is true by \cite{AP08}.  When $m=1$, we could verify the conjecture by computer programming.

\smallskip
Let $f$ be a nilpotent element of $\g=sl_3$. Then $f$ is  minimal or regular.  For $k\in\C$, let $W^k(sl_3,f)$ be the universal affine $W$-vertex algebra associated to $f$ \cite{FF92,KRW03}. Denote by $W_k(sl_3,f)$ the simple quotient of $W^k(sl_3,f)$. 

\smallskip
Let $f$ be a minimal nilpotent element of $\g=sl_3$.  If $k=-3+\frac{2m+1}{2}$, $m\geq 1$, it was shown in \cite{Ara13} that $X_{W_k(sl_3,f)}=\{0\}$ and $W_k(sl_3,f)$ is rational. If $k=-1$, it has been proved in \cite{AKMPP18} that $W_{-1}(sl_3,f)$ is isomorphic to the rank one Heisenberg vertex operator algebra. So $X_{W_{-1}(sl_3,f)}$ is one-dimensional. 
We have the following result for $k=-3+\frac{2}{2m+1}$, $m\geq 1$  (also see Theorem \ref{thm5.1}).
\begin{theorem}\label{main4}
	Let $k=-3+\frac{2}{2m+1}$, $m\geq 1$,  and $f$ a minimal nilpotent element of $sl_3$. Then 
$$
	X_{W_{k}(sl_3,f)}=\left\{\left[\begin{array}{ccc} a & b & 3(\mu^2-a^2)\\
	0 & -2a\mu & c\\
	1 & 0 & a\end{array}\right]| ~\mu, a,b,c\in\C, \ bc=2(a-\mu)(2a+\mu)^2 \right\}.
$$
	In particular,  $\dim X_{W_{k}(sl_3,f)}=3$. Furthermore,  $W_k(sl_3, f)$ is not quasi-lisse. 
\end{theorem}
Let $f$ be a regular nilpotent element of $\g=sl_3$. If $k=-3+\frac{p}{q}$, $p,q\geq 3$, $(p,q)=1$, by \cite{Ar15b},  $X_{W_k(sl_3,f)}=\{0\}$ and $W_k(sl_3,f)$ is rational. Our fifth main result is given as follows (also see Theorem \ref{thm5.2}).
\begin{theorem}\label{main5} 
	Let $k+3=\frac{2}{2m+1}$ or $\frac{2m+1}{2}$,  $m\geq 1$,  and $f$ a regular nilpotent element of $sl_3$. Then
	$$
	X_{W_k(sl_3,f)}=\left\{\left[\begin{array}{ccc} 0 & \frac{3}{4}\mu^2 & \frac{1}{2}\mu^3\\
	2 & 0 & \frac{3}{4}\mu^2\\
	0 & 2 & 0\end{array}\right]| ~\mu\in\C\right\}.
	$$
	In particular, $\dim 	X_{W_k(sl_3,f)}=1$. Furthermore, $W_k(sl_3,f)$ is not quasi-lisse.
\end{theorem}
The rest of this paper is organized as follows. In Section 2, we recollect some concepts and results that will be needed later.  In Section 3, we study singular vectors and maximal ideals of $V^k(sl_n)$ for $k=-n+\frac{n-1}{q}$ with $n\geq 3$, $q\geq 1$, $(n-1,q)=1$. Section 4 is dedicated  to determining the associated varieties of $L_k(sl_3)$ for $k=-3+\frac{2}{2m+1}$, $m\geq 0$. In Section 5, we obtain the associated varieties of  the simple affine $W$-algebras $W_k(sl_3,f)$. 

\subsection*{Acknowledgements} We  appreciate Tomoyuki Arakawa, Dra{\v{z}}en Adamovi{\'c}, and Anne Moreau for invaluable  discussions and comments. In particular, the first author is grateful to Dra{\v{z}}en Adamovi{\'c} for showing  her the  very helpful result of Theorem 8.2 in  his  paper \cite{AMP23} joint with Pierluigi M{\"o}seneder Frajria and Paolo Papi, and to Tomoyuki Arakawara for telling her  the papers  \cite{Fi06} and \cite{BKK24},  and for the discussions.  Results in this paper were presented in part   in conferences ``Representation Theory XVIII" in Dubrovnik, 2023; ``Vertex Algebras , Infinite Dimensional Lie Algebras, and Related Topics'' in Rome, 2023;
``RIMS Joint Research Activity: Vertex Algebras and Related Topics " in Kyoto, 2024; and ``TYMRC-Workshop on Quantum Algebras and Representations" in Shenzhen, 2024. We thank the organizers of these conferences. Finally, 
 we are very grateful to  the referees for invaluable comments and suggestions.

\section{Preliminaries}
\label{sec:affine_vertex_algebras}
\subsection{$C_2$-algebras and associated varieties}
Let $V$ be a vertex algebra over ${\mathbb C}$ \cite{FLM88, FHL93, LL04}.  Recall that a subset $S$ of $V$ is called a \textit{strongly generating set} of $V$ if $V$ is linearly spanned by elements of the form:
$$
\{x^1_{-m_1}\cdots x^s_{-m_s}{\bf 1}| x^i\in S, \  m_i\in{\mathbb Z}_{+}, s\geq 0, \ 1\leq i\leq s\}.
$$
Furthermore, if $S$ is a finite set, then $V$ is called  \textit{finitely strongly generated} by $S$.  

For a vertex algebra $V$,  let $C_2(V)=\on{span}_{\C}\{a_{-2}b\mid a,b\in V\}$ and
$$R_V= V/C_2(V).$$
Then $R_V$ is a Poisson algebra \cite{Zhu96} with relations:
$$\bar a \cdot \bar b = \overline{a_{-1} b}, \qquad \{\bar a, \bar b\}= \overline{a_{0}b},$$
for $a, b \in V$, where $\bar a = a+C_2(V)$. It is easy to see that $V$ is finitely strongly generated 
if and only if $R_V$ is finitely generated.  $V$ is called \textit{$C_2$-cofinite} if $R_V$ is finite-dimensional \cite{Zhu96}.

Let $V$ be a finitely strongly generated vertex algebra. By definition \cite{Ar12a}, the
{\em associated variety} of $V$
is
the reduced scheme $$X_V= {\rm Specm}(R_V).$$

Recall from \cite{Ar12a} that $V$ is called {\em lisse} if $X_V=\{0\}$. We have the following result.
\begin{theorem}[\cite{Ar12a}]
	Let $V$ be a finitely strongly generated vertex algebra. Then $V$ is lisse if and only if $V$ is $C_2$-cofinite.
	\end{theorem}

 Since  $R_V$ is a finitely generated Poisson algebra,  the variety $X_V$ is a union of symplectic leaves \cite{BG03,ArK18}. 
 \begin{defn}[\cite{ArK18}]
 	Let $V$ be a vertex algebra.  $V$ is 
 	called quasi-lisse if $X_V$ has finitely many symplectic leaves.
	\end{defn}
 \begin{theorem}[\cite{ArK18}]\label{ArK18}
 Let $V$ be a finitely strongly generated vertex algebra. If $V$ is quasi-lisse, then $V$ has finitely many irreducible ordinary modules. 
 \end{theorem}

\subsection{Affine vertex algebras}
Let $\g$ be a finite-dimensional simple Lie algebra over ${\mathbb C}$ with the following normalized bilinear form:
$$(\cdot|\cdot)=\displaystyle{\frac{1}{2h^\vee}\times} \text{ Killing form of }\g.$$
Let
$\fg=\g[t,t^{-1}]\oplus  \C K$ be the associated affine Kac-Moody algebra 
 with the commutation relations:
\begin{align*}
[x\otimes t^m,y\otimes t^n]=[x,y]\otimes t^{m+n}+m(x|y)\delta_{m+n,0}K,\quad
[K,\fg]=0,
\end{align*}
where $x,y\in \g$ and $m,n\in \Z$. 
For $x \in \g$ and $m \in \Z$, we will write $x(m)$ for $x
\otimes t^m$.



\label{sec:Universal affine vertex algebras}
For $k \in \C$, set
\begin{align*}
V^k(\g)=U(\fg)\otimes _{U(\g [t]\oplus  \C K)}\C_k,
\end{align*}
where $\C_k$ is the one-dimensional representation of $\g [t]\oplus  \C K$
on which $K$ acts as multiplication by $k$ and $\g[t]$ acts trivially.
By PBW
Theorem, 
we have
\begin{align}
V^k(\g)\cong U(\g\otimes t^{-1}\C[t^{-1}]) = U(t^{-1} \g[t^{-1}]).
\label{eq:PBW}
\end{align}

The space $V^k(\g)$ is naturally graded as follows:
\begin{align*}
V^k(\g) =\bigoplus_{n\in \Z_{\geq 0}}V^k(\g) _{n},
\end{align*}
where the grading is defined by
$$\deg (x^{i_1}(-n_1)\ldots x^{i_r}(-n_r) {\bf 1}) =  \sum_{i=1}^r n_i,
\quad r \geq 0, \; x^{i_j} \in \g,
$$
with ${\bf 1}$ the image of $1\otimes 1$ in $V^k(\g)$.
We have $V^k(\g)_0=\C{\bf 1} $,
and we  may identify $\g$ with $V^k(\g)_1$ via the linear isomorphism
defined by $x\mapsto x(-1){\bf 1} $.

It is well-known that $V^k(\g)$ has a unique vertex algebra structure
such that ${\bf 1}$ is the vacuum vector,
$$x(z)= Y(x\otimes t^{-1},z)  =\sum\limits_{n \in \Z} x(n) z^{-n-1},$$
and
\begin{align*}
[T,x(z)]=\partial_z x(z)
\end{align*}
for $x \in \g$,
where $T$ is the translation operator.
Notice that $x(n)$ acts on $V^k(\g)$ by left multiplication. Then  one can view $x(n)$ as an endomorphism of $V^k(\g)$. 
The vertex algebra
$V^k(\g)$ is called the {\em universal affine vertex algebra}
associated with $\g$ at level $k$ \cite{FZ92,Zhu96,LL04}.

If $k+h^\vee\not=0$, the vertex algebra
$V^k(\g)$ is a vertex operator algebra by the {\em Sugawara construction}.
More specifically, set
$$S=\displaystyle{\frac{1}{2}} \sum_{i=1}^{d}
x_{i}(-1) x^{i}(-1)  {\bf 1},$$
where  $d = \dim \g$, and $\{x_{i}\colon i=1,\ldots,d\}$ is the dual
basis of a basis $\{x^{i}\colon i=1,\ldots,\dim \g\}$ of $\g$
with respect to the bilinear form $(\cdot|\cdot)$.
Then the vector
$\omega=\displaystyle{\frac{S}{k+h^\vee}}$
is a conformal vector of $V^k(\g)$ 
with central charge
$$c(k)=\displaystyle{\frac{k \dim \g}{k+h^\vee}}.$$

Any
graded quotient of $V^k(\g)$
as $\fg$-module
has the structure of a quotient vertex algebra.
In particular,
the  unique simple graded quotient $L_k(\g)$
is a vertex algebra,
and is called the {\em simple affine vertex algebra associated with $\g$ at level $k$. }

\subsection{Zhu's algebra of affine vertex algebras}

For a $\Z_{\geq 0}$-graded vertex algebra $V=\oplus_{n\geq 0}V_n$, let $O(V)$ be the subspace of $V$ linearly spanned by elements
$$
u\circ v=\sum\limits_{j=0}^{\infty}\binom{n}{j} u_{j-2}v,
$$
for $u\in V_n, n\in\Z_{\geq 0}$, $v\in V$. The Zhu's algebra \cite{Zhu96}  
$$
A(V)=V/O(V)
$$
is a unital associative algebra with the multiplication induced from
$$
u*v=\sum\limits_{j=0}^{\infty}\binom{n}{j} u_{j-1}v,
$$
for $u\in V_n, n\in\Z_{\geq 0}$, $v\in V$.   It is known from \cite{FZ92} that the Zhu's algebra of $V^k(\g)$ is the  universal enveloping algebra $U(\g)$ of $\g$.  Let ${\mathcal I}_k$ be the maximal ideal of $V^k(\g)$. Then
$$
L_k(\g)=V^k(\g)/{\mathcal I}_k.
$$
Denote by $I_k$ the image of ${\mathcal I}_k$ in  $A(V^k(\g))$, then
$$
A(L_k(\g))=U(\g)/I_k.
$$

Let $\h$ be the Cartan subalgebra of $\g$, and $\g=\frak{n}^++\h+\frak{n}^-$ be a fixed triangular decomposition of $\g$. Set
$$
U(\g)^{\h}=\{u\in U(\g) | [h, u]=0, \ {\rm for~all} \ h\in \h\}
$$
and let
$$
\frak{p}:  U(\g)^{\h}\to U(\h)
$$
be the Harish-Chandra projection map,  which is the restriction of the projection map: $U(\g)=U(\h)\oplus (\n^-U(\g)\oplus U(\g)\n^+)\to U(\h)$ to $U(\g)^{\h}$. It is known that $\frak{p}$ is an algebra homomorphism \cite{Hu72}.  For a two-sided ideal $I$ of $U(\g)$, the {\em characteristic variety of} $I$ is defined as \cite{Jo77}
$$
{\mathcal V}(I)=\{\lambda\in\h^*|p(\lambda)=0 \ {\rm for~all} \ p\in\frak{p}(I\cap U(\g)^{\h})\}.
$$
For $\lambda\in\h^*$, let $L(\lambda)$ be  the irreducible module of $V^k(\g)$ induced from the simple highest weight $\g$-module $L(\bar{\lambda})$. Then by \cite{Zhu96}, $L(\lambda)$ is a module of $L_k(\g)$ if and only if $I_k$ acts on $L(\bar{\lambda})$ as zero. We have the following result from \cite{Ar16}.
\begin{prop}[\cite{Ar16}]
	$L(\lambda)$ is an $L_k(\g)$-module if and only if $\bar{\lambda}\in{\mathcal V}(I_k)$.
\end{prop}
\subsection{Associate graded vertex Poisson algebras of affine vertex algebras}
It is known by Li \cite{Li05}
that any vertex algebra $V$ admits a canonical filtration $F^\bullet V$,
called the {\em Li filtration} of $V$.

\vskip 0.2cm
For a quotient $V$ of $V^k(\g)$,
$F^\bullet V$ is described as follows.
The subspace
 $F^p V$ is spanned by the elements
$$
y_{1}(-n_1-1)\cdots y_{r}(-n_r-1)\mathbf{1}
$$
with $y_{i} \in \g$,
$n_i\in\Z_{\geq 0}$, $ n_1+\cdots +n_r\geq p$.
We have
\begin{align*}\nonumber
& V=F^0V\supset F^1V\supset\cdots, \quad \bigcap_{p}F^pV=0,\\ 
& TF^pV\subset F^{p+1}V,\\\nonumber
& a_{n}F^{q}V\subset F^{p+q-n-1}V ~\text{for}~a\in F^{p}V, \ n\in\Z,\\\nonumber
& a_{n}F^{q}V\subset F^{p+q-n}V ~\text{for}~a\in F^{p}V, \ n\geq 0.
\end{align*}
Here we set $F^pV=V$ for $p<0$.

Let ${\rm gr}^FV=\bigoplus_p F^pV/F^{p+1}V$ be the associated graded vector space.
The space ${\rm gr}^FV$ is a vertex Poisson algebra by
\begin{align*}
& \sigma_{p}(a)\sigma_{q}(b)=\sigma_{p+q}(a_{-1}b),\\
& T\sigma_{p}(a)=\sigma_{p+1}(Ta),\\
& \sigma_{p}(a)_{n}\sigma_{q}(b)=\sigma_{p+q- n}(a_{n}b)
\end{align*}
for $a,b\in V$,
$n\geq 0$,
where $\sigma_p \colon F^p(V)\rightarrow F^pV/F^{p+1}V$ is the principal symbol map.
In particular, the space
$\on{gr}^F V$ is a $\g[t]$-module by the correspondence
\begin{align}
\g[t]\ni x(n)\longmapsto \sigma_0(x)_{}\in \End(\on{gr}^F V)
\label{eq:g[t]-action-in-gr}
\end{align}
for $x\in \g$, $n\geq 0$.

The filtration  $F^\bullet V$
is compatible with the grading:
$F^pV=\bigoplus\limits_{n\in \Z_{\geq 0}}F^pV_{n}$,
where
 $F^p V_n := V_n \cap F^p V$.

\subsection{Zhu's $C_2$-algebras and associated varieties
of affine vertex algebras}
\label{sub:Zhu affine vertex algebras}

Let $F^{\bullet}V$ be given as above. 
We have \cite[Lemma 2.9]{Li05}
$$F^p V = {\rm span}_\C\{ a_{-i-1} b \colon a \in V, i \geq 1, b \in F^{p-i} V \}$$
for all $p \geq 1$. In particular, 
$$F^1 V = C_2(V),\qquad R_V = V/C_2(V) = F^0 V / F^1 V \subset {\rm gr}^FV.$$
%

It is easy to see  that
$$F^1 V^k(\g) = C_2(V^k(\g)) = t^{-2} \g[t^{-1}] V^k(\g).$$

The following map defines an isomorphism of Poisson algebras
\begin{align*}
\begin{array}[t]{rcl}
\C[\g^*] \cong S(\g) & \longrightarrow & R_{V^k(\g)},  \\[0.2em]
\g \ni x  & \longmapsto  & x (-1) {\bf 1}
+ t^{-2} \g[t^{-1}] V^k(\g).
\end{array}
\end{align*}
So $$R_{V^k(\g)}  \cong \C[\g^*]$$ and 
$$
X_{V^k(\g)}\cong \g^*.$$
More generally, let $V$ be a quotient of $V^k(\g)$
by an ideal $N$, then
we have
\begin{align}
R_{V}\cong \C[\g^*]/\overline{N}
\end{align}
as Poisson algebras,
where $\overline{N}$ is the image of $N$ in $R_{V^k(\g)}=\C[\g^*]$. 
Then $X_V$ is just the zero locus of $\overline{N}$ in
$\g^*$.
It is a closed $G$-invariant conic subset of $\g^*$. Identifying $\g^*$ with $\g$ via the bilinear form $(\cdot|\cdot)$,
one may view $X_V$ as a subvariety of $\g$. 
\subsection{Sheets and Jordan classes}

In this subsection, we recall Jordan classes and sheets of a semisimple Lie algebras following \cite{TY05}. 
Let $\frak{g}$ be as above. Denote by ${\mathcal S}$ (resp. ${\mathcal N}$) the set of semisimple (resp. nilpotent) elements of ${\frak g}$.  
For $x\in\frak{g}$, set $\mathfrak{g}^x=\{y\in\frak{g}|[x, y]=0\}$, and denote by $x_s$ and  $x_n$ its semisimple and nilpotent components, respectively.  
For $x,y\in\frak{g}$, we say that $x$ and $y$ are $G$-Jordan equivalent if there exists $\alpha\in G$ such that 
$$
\frak{g}^{y_s}=\frak{g}^{\alpha(x_s)}=\alpha(\frak{g}^{x_s}), \qquad y_n=\alpha(x_n).
$$
This defines an equivalence relation on $\frak{g}$. The equivalence class of $x$, denoted by $J_G(x)$, is called the {\it Jordan class} of $x$ in $\frak{g}$.  A Jordan class is clearly a $G$-stable set. 

For any Lie subalgebra $\frak{u}\subset\frak{g}$, set
$$
\frak{z}(\frak{u})=\{x\in\frak{u}|[x, y]=0, \ y\in\frak{u}\}.
$$
Denote by $\frak{u}^{reg}$   the set of elements $y\in\frak{u}$ such that the dimension of $\g^y$ is minimal. We have the following lemma.

\begin{lem}
	Let $x\in\frak{g}$.  
	\begin{enumerate}
		\item We have
		$$
		J_G(x)=G(\frak{z}(\frak{g}^{x_s})^{reg}+x_n).
		$$
		\item  $J_G(x)$ is irreducible in $\frak{g}$. 
		\item  $J_G(x)$ is locally closed in $\frak{g}$, so it is a subvariety of $\frak{g}$. 
	\end{enumerate}
\end{lem}
To a Jordan class $J$, let  $x\in J$,  and  $\frak{l}=\frak{g}^{x_s}$. Then $\frak{l}$  is a Levi subalgebra of $\frak{g}$.  Let ${\mathbb O}_{\frak l}$ be the nilpotent orbit in $\frak{l}$ of $x_n$.  The pair $(\frak{l}, {\mathbb O}_{\frak{l}})$ does not depend on $x\in J$ up to $G$-conjugacy, and there is a one-to-one correspondence between the set of pairs $(\frak{l}, {\mathbb O}_{\frak{l}})$ and the set of Jordan classes. 

\vskip 0.3cm
Let $n\in{\mathbb N}$. Denote
$$
\ {\frak g}^{(n)}=\{x\in\frak{g}|\dim{\frak g}^{x}=n\}.
$$
\begin{defn} For $n\in\N$, an irreducible component of $\frak{g}^{(n)}$ is called a {\it sheet} of $\frak{g}$.  
\end{defn}
A sheet of $\g$ is a disjoint union of Jordan classes. So a sheet ${\mathbb S}$ contains a unique dense open Jordan class $J$. The datum $(\frak{l}, {\mathbb O}_{\frak l})$ of the Jordan class $J$ is called the {\it datum} of ${\mathbb S}$. Then 
$$
\overline{\mathbb S}=\overline{J}, \qquad {\mathbb S}=(\overline{J})^{reg}.
$$
Let ${\mathbb S}$ be a sheet with datum $(\frak{l}, {\mathbb O}_{\frak l})$, then the induced nilpotent orbit ${\rm Ind}_{\frak l}^{\g}({\mathbb O}_{\frak l})$ of $\g$ from ${\mathbb O}_{\frak l}$ in ${\frak l}$ is the unique orbit contained in ${\mathbb S}$.  The rank of ${\mathbb S}$ with datum $(\frak{l}, {\mathbb O}_{\frak l})$ is by definition
\begin{equation}\label{sheet1}
r({\mathbb S}):=\dim {\mathbb S}-\dim ({\rm Ind}_{\frak l}^{\g}({\mathbb O}_{\frak l}))=\dim\frak{z}({\frak l}).
\end{equation}
A sheet ${\mathbb S}$ with datum $(\frak{l}, {\mathbb O}_{\frak l})$ is called a {\it Dixmier sheet} if ${\mathbb O}_{\frak l}=\{0\}$.

\subsection{Affine $W$-algebras}
Let $f$ be a nilpotent element of $\g$. By the Jacobson-Morozov Theorem,
there is an ${sl}_2$-triple $(e,h,f)$ of $\g$.
Recall that the Slodowy slice $\mathscr{S}_f$ is the affine space $f+\g^{e}$.
It has a natural Poisson structure induced from that of $\g^*$ \cite{GG02}.

The embedding ${\rm span}_\C\{e,h,f\}
\cong {sl}_2 \hookrightarrow \g$
exponentiates to a homomorphism
$SL_2 \to G$. By restriction to the one-dimensional
torus consisting of diagonal matrices, we obtain
a one-parameter subgroup $\rho \colon \C^* \to G$.
For $t\in\C^*$ and $x\in\g$, set
\begin{align*} \label{eq:rho}
\tilde{\rho}(t)x := t^{2} \rho(t)(x).
\end{align*}
We have
$\tilde{\rho}(t)f=f$, and
the $\C^*$-action of $\tilde{\rho}$ stabilizes $\mathscr{S}_f$.
Moreover, it contracts to $f$ on $\mathscr{S}_f$, that is, for all $x\in\g^{e}$,
$$\lim_{t\to 0} \tilde{\rho}(t)(f+x)=f.$$

The following proposition is well-known.
\begin{prop}[{\cite{S80,P02,CM16}}]
	\label{pro:Slodowy}
	The morphism
	$$\theta_f \colon G \times \mathscr{S}_f \longrightarrow \g,
	\quad (g,x) \longmapsto g.x$$
	is smooth onto a dense open subset of $\g^*$.
\end{prop}

Let $W^k(\g,f)$ be the {\it affine $W$-algebra} associated with
a nilpotent element $f$ of $\g$
defined by the generalized quantized Drinfeld-Sokolov reduction:
$$W^k(\g,f)=H^{0}_{DS,f}(V^k(\g)).$$
Here $H^{\bullet}_{DS,f}(M)$ denotes the BRST
cohomology of the  generalized quantized Drinfeld-Sokolov reduction
associated to $f \in \mathcal{N}(\g)$ with coefficients in
a $V^k(\g)$-module $M$ \cite{KRW03}.
Recall that we have a natural isomorphism
$R_{W^k(\g,f)}\cong \C[\mathscr{S}_{f}]$ of Poisson algebras  \cite{DSK06,Ar15a}, so that
\begin{align*}
X_{W^k(\g,f)}= \mathscr{S}_{f}.
\end{align*}
We write $W_k(\g,f)$ for the unique simple (graded) quotient of
$W^k(\g,f)$. Then $X_{W_k(\g,f)}$
is a $\C^*$-invariant Poisson
subvariety of the Slodowy slice $\mathscr{S}_f$.

Let $\mathcal{O}_k$ be the category $\mathcal{O}$ of
$\fg$ at level $k$.
We have a functor
\begin{align*}
\mathcal{O}_k\longrightarrow W^k(\g,f)\on{-Mod}
,\quad M\longmapsto
H^0_{DS,f}(M),
\end{align*}
where
$W^k(\g,f)\on{-Mod}$ denotes the category
of $W^k(\g,f)$-modules.

The full subcategory of $\mathcal{O}_k$ consisting of
objects $M$ on which $\g$ acts  locally finitely is denoted by $\on{KL}_k$.
Note that both $V^k(\g)$ and $L_k(\g)$ are objects of $\on{KL}_k$.
\begin{theorem}[\cite{Ar16}]\label{Ar16}
	 Let 
	 $f_{\theta}$ be a root vector of the highest root $\theta$ of $\g$. Then  
	 	\begin{enumerate}
	 	\item $H_{DS,f_{\theta}}^i(M)=0$ for all $i\ne 0$,  and $M\in
	 	{\mathcal O}_k$.
	 	In particular, the functor
	 	\begin{align*}
	 		\mathcal{O}_k\longrightarrow W^k(\g,f_{\theta})\on{-Mod}
	 		,\quad M\longmapsto
	 		H^0_{DS,f_{\theta}}(M),
	 	\end{align*}
 		 is exact.
	 	\item 
	 	$H^0_{DS,f_{\theta}}(L(\lambda))\neq 0$ if and only if $\lambda(\al_0^{\vee})\notin\Z_{\geq 0}$, where $\al_0^{\vee}=-\theta^{\vee}+K$.  If this is the case, $H^0_{DS,f_{\theta}}(L(\lambda))$ is a simple $W^k(\g,f_{\theta})$-module.
	 	\end{enumerate}
	\end{theorem}
\begin{theorem}[{\cite{Ar15a}}]
	\label{Th:W-algebra-variety}
	Let $k, f$ be arbitrary. Then
	\begin{enumerate}
		\item $H_{DS,f}^i(M)=0$ for all $i\ne 0$, and $M\in
		\on{KL}_k$.
		In particular, the functor
		\begin{align*}
			\on{KL}_k\longrightarrow W^k(\g,f)\on{-Mod}
			,\quad M\longmapsto
			H^0_{DS,f}(M),
		\end{align*}
		is exact.
		\item
		For any quotient $V$ of $V^k(\g)$,
		\begin{align*}
		X_{H^{0}_{DS,f}(V)}=X_{V}\cap  \mathscr{S}_{f}.
		\end{align*}
		In particular,
		$H_{DS,f}^{0}(V)
		\ne 0$ if and only if
		$\overline{G.f}\subset X_V$.
		\label{item:intersection}
	\end{enumerate}
\end{theorem}

\section{Maximal ideals of $V^k(sl_3)$ for $k+n=\frac{n-1}{q}$,  $n\geq 3$, $q\geq 2$, $(q,n-1)=1$}

\vskip 0.2cm
 Let $\g$, $\fg$  and $V^k(\g)$ be the same as  in Section 2.   We first recall the following result from \cite{GK07}.
\begin{theorem}[\cite{GK07}]\label{GK07}  For $k\in\C$, 
	$V^k(\g)$
	is not simple if and only if
	\begin{equation}	\label{eq:Gorelik-Kac}
	r^{\vee}(k+h^{\vee}) \in {\mathbb Q}_{\geq 0}\backslash \left\{\frac{1}{m}\mid m\in {\mathbb Z}_{\geq 1}\right\},
	\end{equation}
	where
	$r^{\vee}$ is the lacing number of $\g$.
\end{theorem}

By Theorem \ref{GK07},  $V^k(sl_n)$ is not simple if and only if $k=-n+\frac{p}{q}$, where $p\in\Z_{\geq 2}$, $q\in\Z_{\geq 1}$, and $(p,q)=1$. 

\vskip 0.2cm
In this section, we will prove that, for  $k=-3+\frac{2}{2m+1}$, $m\in\Z_{\geq 1}$, the maximal ideal of $V^k(sl_3)$ is generated by two singular vectors of conformal weight $3(2m+1)$, and for $n\geq 4$, $k=-n+\frac{n-1}{q}$, $q\geq 2$, $(q,n-1)=1$, the maximal ideal of $V^k(sl_n)$ is generated by a singular vector of conformal weight $2q$. We reach our results by using   the Kashiwara-Tanisaki characters \cite{KT00},  category equivalences established  in \cite{Fi06}, the character formulas for the $k=-1$ case given in \cite {KW18, BKK24},  and the  results  established in \cite{AP08, AKMPP20, AM17}.

\subsection{Kashiwara-Tanisaki Characters}

In this subsection, we will introduce the Kashiwara-Tanisaki characters following \cite{KT00}.  Let $\Delta$, $\widehat{\Delta}$ be the root systems of $\g$ and $\fg$ with respect to $\h$ and $\fh$, respectively. Denote by $\delta$ the positive
imaginary root such that any imaginary root is an integral multiple of $\delta$. 
Let $W$ and $\wW$ be the Weyl groups of $\g$ and $\fg$,  and $Q$ and $\wQ$ the root lattices of $\g$ and $\fg$, respectively. Let $\{\al_1,\dots,\al_l\}$ be the simple root system of $\g$, and $\{\al_0,\al_1,\dots,\al_l\}$  the simple root system of $\fg$.  For a real root $\al\in\wD$, set $\al^{\vee}=2\al/(\al|\al)$.   Let $\rho\in{\frak h}^*$ and $\wr\in\fh^*$ be such that 
$$
(\rho|\al_i^{\vee})=1,   \   \  (\wr|\al_0^{\vee})=(\wr|\al_i^{\vee})=1, \ 1\leq i\leq l.
$$
Recall that the {\it twisted (or shifted) action} of $\wW$ on $\fh^*$ is defined as follows
$$
w\circ \lambda=w(\lambda+\wr)-\wr,
$$
where $w\in \wW$ and $\lambda\in\fh^*$.  Notice that if $w\in W$, $w\circ \lambda=w(\lambda+\rho)-\rho$.

 For $\lambda\in\fh^*$, set 
$$
\wD(\lambda)=\{\alpha\in\wD^{re}|(\lambda+\wr|\al^{\vee})\in\Z\}, 
$$
$$
\wD_0(\lambda)=\{\alpha\in\wD^{re}|(\lambda+\wr|\al^{\vee})=0\}.
$$
Notice that $\wD(\lambda)$ and $\wD_0(\lambda)$ are subsystems of $\wD^{re}$. Denote the set  of positive roots, the set of negative roots, the set of simple roots  and the Weyl group for  $\wD(\lambda)$ by
$\wD^+(\lambda)$,  $\wD^-(\lambda)$, $\widehat{\Pi}(\lambda)$ and $\wW(\lambda)$, respectively. Denote those for $\wD_0(\lambda)$ by  $\wD^{+}_0(\lambda)$, $\wD^{-}_0(\lambda)$, $\widehat{\Pi}_0(\lambda)$ and $\wW_0(\lambda)$.

For a real root $\al\in\wD$,  denote by $s_{\al}\in\wW$ the corresponding reflection. Then  $\widehat{\Pi}(\lambda)$ is  the set of $\al\in\wD^+(\lambda)$ such that $s_{\al}(\wD^+(\lambda)\backslash \{\al\})=\wD^+(\lambda)\backslash \{\al\}$, and  $(\wW(\lambda), S(\lambda))$ is a Coxeter group, where   $S(\lambda)=\{s_{\al}: \al\in\widehat{\Pi}(\lambda)\}$ \cite{KW89}, \cite{KT00}. 

\vskip 0.2cm
For $w\in \wW(\lambda)$, denote by $l_{\lambda}(w)$ the length of $w$. Denote the Bruhat ordering of $\wW(\lambda)$ by $\geq_{\lambda}$. For $y,w\in \wW(\lambda)$, denote by $P^{\lambda}_{y,w}(q)\in\Z[q]$ the associated Kazhdan-Lusztig polynomial  \cite{KL79}, and by $Q^{\lambda}_{y,w}(q)\in\Z[q]$ the inverse Kazhdan-Lusztig polynomial defined by 
$$
\sum\limits_{x\leq_{\lambda}y\leq_{\lambda}z}(-1)^{l_{\lambda}(y)-l_{\lambda}(x)}Q^{\lambda}_{x,y}(q)P^{\lambda}_{y,z}(q)=\delta_{x,z},
$$
for any $x,z\in \wW(\lambda)$.  Set
$$
{\mathcal C}=\{\lambda\in\fh^*| (\delta|\lambda+\wr)\neq 0 \}, 
$$
$$
{\mathcal C}^+=\{\lambda\in{\mathcal C}|(\lambda+\wr|\al^{\vee})\geq 0, \  {\rm for \ any} \ \al\in\wD^{+}(\lambda)\}, 
$$
and 
$$
{\mathcal C}^-=\{\lambda\in{\mathcal C}|(\lambda+\wr|\al^{\vee})\leq 0, \  {\rm for \ any} \ \al\in\wD^{+}(\lambda)\}.
$$
Let $\lambda\in{\mathcal C}$. Then $\wW_0(\lambda)$ is a finite group \cite{KT00}. 

For $\lambda\in\fh^*$, let $M(\lambda)$ (resp. $L(\lambda)$) be the Verma module (resp. simple module) of $\fg$  with highest weight $\lambda$. 
 The following lemma comes from \cite{KK79},  \cite{DGK82}, \cite{Ku87},  \cite{KW89}.
\begin{lem}\label{lweyl3} Let $\lambda\in{\mathcal C}^+$, and $w\in\wW(\lambda)$. If $L(\mu)$ is a subquotient of $M(w\circ \lambda)$, then there exists $y\in \wW(\lambda)$ such that $\mu=yw\circ \lambda$, and $w\leq_{\lambda}yw$.
\end{lem}	
We have the following character formulas from \cite{KT00}.
\begin{theorem}\cite{KT00}\label{KT}
	\begin{enumerate}
		\item 	Let $\lambda\in{\mathcal C}^+$,  then for any $w\in\wW(\lambda)$ which is the longest element of $w\wW_0(\lambda)$, 
		$$
		{\rm ch}(L(w\circ \lambda))=\sum\limits_{w\leq_{\lambda}y\in \wW(\lambda)}(-1)^{l_{\lambda}(y)-l_{\lambda}(w)}Q^{\lambda}_{w,y}(1){\rm ch}(M(y\circ \lambda)).
		$$
		\item  	Let $\lambda\in{\mathcal C}^-$,  then for any $w\in\wW(\lambda)$ which is the shortest element of $w\wW_0(\lambda)$, 
		$$
		{\rm ch}(L(w\circ \lambda))=\sum\limits_{w\geq_{\lambda}y\in \wW(\lambda)}(-1)^{l_{\lambda}(w)-l_{\lambda}(y)}P^{\lambda}_{y,w}(1){\rm ch}(M(y\circ \lambda)).
		$$
	\end{enumerate}  
\end{theorem}
We have the following lemma from \cite{KT98}, \cite{Fi06}.
\begin{lem}\label{hom}
	Let $\lambda\in{\cal C}^+$ or ${\cal C}^-$, and  $w, w'\in\wW(\lambda)$ be such that $w'\cdot \lambda\leq w\cdot \lambda$. Then
	\begin{equation}\label{ehom}
		{\rm dim}{\rm Hom}_{\mathbb C}(M(w'\cdot \lambda), M(w\cdot \lambda))=1.
		\end{equation}
\end{lem}

The following theorem comes from \cite{Fi06}.
\begin{theorem}\label{Fi06}
	Let ${\frak h}\subset{\frak b}\subset{\g}$ and ${\frak h}'\subset{\frak b}'\subset{\g}'$ be symmetrizable Kac-Moody Lie algebras with Cartan and Borel subalgebras, respectively. Let  $\Lambda\in\frak{h}^*/\sim$ and $\Lambda'\in(\frak{h}')^*/\sim'$ be two equivalent classes outside the critical hyperplanes and let ${\cal O}_{\Lambda}$ and ${\cal O}'_{\Lambda'}$ be the corresponding indecomposable blocks, where $\sim$ and $\sim'$ are the usual relation in $\h^*$ and $(\h')^*$, respectively (see \cite{KK79}). Suppose the following conditions hold:
	\begin{enumerate}
	\item \ There exist $\lambda\in\Lambda$ and $\lambda'\in\Lambda'$ which are either  both dominant or antidominant; 
	
	\item $(W(\lambda), S(\lambda))\cong (W(\lambda'), S(\lambda'))$ are isomorphic as Coxeter groups;
	
	\item $W_0(\lambda)\cong W_{0}(\lambda')$ under the same isomorphism and both sets are finite.
	\end{enumerate}
Then there exists an equivalence of categories
$$
{\cal O}_{\Lambda}\cong {\cal O'}_{\Lambda'}.
$$
\end{theorem}
\begin{remark}
	As an application of Theorem \ref{Fi06}, the author in \cite{Fi06} proved  Theorem  \ref{KT} for symmetrizable  Kac-Moody Lie algbera $\g$,  and regular $\Lambda\in{\frak h}^*/\sim$ satisfying the condition in Theroem \ref{Fi06} and that $(W(\Lambda), S(\Lambda))$ is finite or affine type,   and dominant or antidominant $\lambda\in\Lambda$.
\end{remark}

\subsection{Maximal ideals of $V^k(sl_n)$ for $k=-n+\frac{n-1}{q}$($q\geq 2$, $(q,n-1)=1)$}
In this subsection, we assume that 
$\g=sl_n(\C)$, and $k=-n+\frac{n-1}{q}$, $q\geq 1$, $(q, n-1)=1$. It is obvious that
$$
\widehat{\Pi}(k\Lambda_0)=\{\beta_0=q\delta-\theta, \al_1, \cdots, \al_{n-1}\}.
$$
Since $(k\Lambda_0+\wr|\beta_0)=0$, it follows that
$$
\widehat{\Pi}_0(k\Lambda_0)=\{\beta_0\},
$$
and
$$
\wW_0(k\Lambda_0)=\{e, s_{\beta_0}\}.
$$
In particular, if $q=1$, then $k=-1$ and $\beta_0=\al_0$. It follows that
 $$\wW(-\Lambda_0)=\wW, \ \wW_0(-\Lambda_0)=\{e, s_{\al_0}\}.$$
   Let $S=\{s_{\al_0}, s_{\al_1}, \cdots, s_{\al_{n-1}}\}$, then $(\wW, S)$ is a Coxeter group. 
We have the following lemma.
\begin{lem}\label{lweyl1} 
	For any $k=-n+\frac{n-1}{q}$, $q\geq 1$, $n\geq 3$, we have
	$$
	(\wW(k\Lambda_0), S(k\Lambda_0))\cong (\wW, S)
	$$
	as Coxeter groups by $\sigma: (\wW, S)\to (\wW(k\Lambda_0), S(k\Lambda_0))
	$
	such that
	$$\sigma s_{\al_0}=s_{\beta_0}, \sigma s_{\al_i}=s_{\al_i}, \ i=1,2,\cdots,n-1.	
	$$
	In particular, $\sigma(W_0(-\Lambda_0))=W_0(k\Lambda_0)$.
\end{lem}
\begin{proof} Recall that  the defining relations of $(\wW, S)$ are \cite{K90}
	$$
	s_{\al_i}^2=1, \ i=0,1,\cdots,n-1,  \  (s_{\al_i}s_{\al_j})^{m_{ij}}=1 \  i,j=0,1,\cdots, n-1,
	$$
such that 
$$
m_{01}=m_{10}=m_{0,n-1}=m_{n-1,0}=3, \  m_{0j}=m_{j0}=2, \ j=2,\cdots,n-2,
$$
$$
m_{i,i+1}=m_{i+1,i}=3, \ m_{ij}=m_{ji}=2, \ 1\leq i,j\leq n-1, \ j\neq i+1. 
$$
Notice that $S=\{s_{\al_0}, s_{\al_1}, \cdots, s_{\al_{n-1}}\}$, $S(k\Lambda_0)=\{s_{\beta_0}, s_{\al_1}, \cdots, s_{\al_{n-1}}\}$. 
 Thus it is enough to prove that 
\begin{equation}\label{eweyl1}
	(s_{\beta_0}s_{\al_1})^3=(s_{\al_1}s_{\beta_0})^3=(s_{\beta_0}s_{\al_{n-1}})^3=(s_{\al_{n-1}}s_{\beta_0})^3=1,
\end{equation}
and
\begin{equation}\label{eweyl2}
(s_{\beta_0}s_{\al_j})^2=(s_{\al_j}s_{\beta_0})^2=1,  \ j=2,\cdots, n-2.
\end{equation}
We only check 	$(s_{\beta_0}s_{\al_1})^3=1$. The computation  for
 the other cases in (\ref{eweyl1})-(\ref{eweyl2}) is similar. Notice that
 $$
 s_{\beta_0}s_{\al_1}s_{\beta_0}=s_{s_{\beta_0}(\al_1)}=s_{q\delta-\theta+\al_1},
 $$
 and
 $$
 s_{\al_1}s_{q\delta-\theta+\al_1}s_{\al_1}=s_{s_{\al_1}(q\delta-\theta+\al_1)}=s_{q\delta-\theta}=s_{\beta_0}.
 $$
 Then we have
 $$
 (s_{\beta_0}s_{\al_1})^3=s_{\beta_0}s_{\al_1}s_{\beta_0}s_{\al_1}s_{\beta_0}s_{\al_1}=s_{\beta_0}s_{\al_1}s_{q\delta-\theta+\al_1}s_{\al_1}=s_{\beta_0}s_{\beta_0}=1.
 $$
	\end{proof}
For $\gamma\in \h^*$, let $t_{\gamma}\in{\rm End}\fh^*$ be defined by \cite{K90}
$$
t_{\gamma}(\lambda)=\lambda+(\lambda|\delta)\gamma-((\lambda|\gamma)+\frac{1}{2}(\lambda|\delta)(\gamma|\gamma))\delta, 
$$
for $\lambda\in\fh^*$. 

\vskip 0.2cm
Recall from \cite{K90}, for any $x\in\wW$, there exists $y\in W$ and $\gamma\in Q$ such that $x=yt_{\gamma}$.  We have the following lemma.
\begin{lem}\label{lweyl2} 
	For $\gamma\in Q$,  $\sigma(t_{\gamma})=t_{q\gamma}$. 
	\end{lem}
\begin{proof} We first consider $\gamma=\theta$. Notice that for $\lambda\in\fh^*$, 
		$$
	s_{\beta_0}s_{\theta}(\lambda)=s_{\beta_0}(\lambda-(\lambda|\theta)\theta)=\lambda+(\lambda|\delta)q\theta-((\lambda|\theta)q+((\lambda|\delta)q^2)\delta=t_{q\theta}(\lambda).
	$$
	This deduces that
	$$
	s_{\beta_0}s_{\theta}=t_{q\theta}.
	$$
In particular, if $q=1$,  $s_{\al_0}s_{\theta}=t_{\theta}$ (see also \cite{K90}). So by Lemma \ref{lweyl1}, 
\begin{equation}\label{eweyl3}
\sigma(t_{\theta})=t_{q\theta}.
\end{equation}
Notice that $$s_{\al_1}t_{\theta}s_{\al_1}=t_{s_{\al_1}(\theta)}=t_{\theta-\al_1}$$
 and 
 $$s_{\al_1}t_{q\theta}s_{\al_1}=t_{s_{\al_1}(q\theta)}=t_{q(\theta-\al_1)}.$$
 Thus
 \begin{equation}\label{eweyl4}
\sigma(t_{\theta-\al_1})=t_{q(\theta-\al_1)}.
\end{equation}
Considering $s_{\al_2}t_{\theta-\al_1}s_{\al_2}$ and $s_{\al_2}t_{q(\theta-\al_1)}s_{\al_2}$, we  deduce that
\begin{equation}\label{eweyl5}
\sigma(t_{\theta-\al_1-\al_2})=t_{q(\theta-\al_1-\al_2)}.
\end{equation}
Similarly, we have for $4\leq j\leq n-1$, 
\begin{equation}\label{eweyl6}
\sigma(t_{\al_j+\cdots \al_{n-1}})=t_{q(\al_j+\cdots \al_{n-1})}.
\end{equation}
Notice that for $\al, \beta\in Q$, $t_{\al}t_{\beta}=t_{\al+\beta}$. Then the lemma follows from  (\ref{eweyl3})-(\ref{eweyl6}).
\end{proof}

Let $\Lambda_i\in\fh^*$ be such that
$$
\Lambda_i(\al_j^{\vee})=\delta_{ij}, \ i,j=0,1,\cdots,n-1,
$$
We have the following character formulas from  \cite{KW18}, \cite{BKK24}. 
\begin{theorem}[\cite{KW18,BKK24}\label{Kaccharacter}]
Let $n\geq 3$, and let $L(\Lambda)$ be an irreducible	$\widehat{sl}_n$-module  at  level $k=-1$ with highest weight 
$$
\Lambda=-(1+s)\Lambda_0+s\Lambda_{n-1} \ ({\rm resp.} \ \Lambda=-(1+s)\Lambda_0+s\Lambda_{1}), \ s\in\Z_{\geq 0}.
$$
Then the character of $L(\Lambda)$ is given  by the following formula:
$$
\widehat{R}{\rm ch}L(\Lambda)=
\sum\limits_{w\in W}\sum\limits_{\tiny{\begin{gathered}\gamma\in Q \\(\gamma|\bar{\Lambda}_{n-1}(resp. \bar{\Lambda}_1))\geq 0\end{gathered}}}(-1)^{l(w)}e^{wt_{\gamma}({\Lambda+\wr})-\wr},
$$
where
$$
\widehat{R}=\prod_{\al\in\wD_+}(1-e^{-\al})^{{\rm dim}\fg_{\al}}.
$$
In particular, if $\Lambda=-\Lambda_0$, 
\begin{equation}\label{emaximal26}
\widehat{R}{\rm ch}L(-\Lambda_0)=
\sum\limits_{w\in W}\sum\limits_{\tiny{\begin{gathered}\gamma\in Q \\(\gamma|\bar{\Lambda}_{n-1})\geq 0, (\gamma|\bar{\Lambda}_1))\geq 0\end{gathered}}}(-1)^{l(w)}e^{wt_{\gamma}(-\Lambda_0+\wr)-\wr}.
\end{equation}
\end{theorem}	

\vskip 0.2cm
Recall  the following results from \cite{AP08},  \cite{AKMPP20}, and \cite{AM17}.
\begin{theorem}\label{AP08}
	\begin{enumerate}
		\item  \cite{AP08}	Let
		$$
		\begin{array}{ll}
			u^1=& [-e_{\al_1}(-1)^2e_{\al_2}(-1)+e_{\al_1}(-1)e_{\al_1+\al_2}(-1)h_2(-1)+e_{\al_1+\al_2}(-1)^2f_{\al_2}(-1)]{\bf 1},\\\\
			u^2=&[e_{\al_1}(-1)e_{\al_2}(-1)^2+e_{\al_2}(-1)e_{\al_1+\al_2}(-1)h_1(-1)\\
			&-2e_{\alpha_2}(-1)e_{\al_1+\al_2}(-2)
			-e_{\al_1+\al_2}(-1)^2f_{\al_1}(-1)]{\bf 1}.
		\end{array}
		$$
		Then $u^1, u^2$ are  singular vectors of $V^{-1}(sl_3)$. 
		
		\item \cite{AP08} Let
		$$
		u=e_{\theta}(-1)e_{\al_2+\cdots+ \al_{n-2}}(-1){\bf 1}-e_{\al_1+\cdots +\al_{n-2}}(-1)e_{\al_2+\cdots +\al_{n-1}}(-1){\bf 1}.
		$$
		Then $u$ is a singular vector of $V^{-1}(sl_n)$, $n\geq 4$. 
		
		\item \cite{AKMPP20} The ideal of $V^{-1}(sl_3)$ generated by $u^1$ and $u^2$ is the maximal ideal of $V^{-1}(sl_3)$.

		\item \cite{AM17} The ideal  generated by $u$ in (2) is the maximal ideal of $V^{-1}(sl_n)$, $n\geq 4$.
	\end{enumerate}
	
\end{theorem}

We are now in a position to state the main result of this section.
\begin{theorem}\label{tmain3.1}  Let $k=-n+\frac{n-1}{q}$, $n\geq 3$, $q\geq 2$, $(n-1,q)=1$. Then 	
\begin{equation}\label{emaximal3}
\widehat{R}{\rm ch}L_k(sl_n)=
\sum\limits_{w\in W}\sum\limits_{\tiny{\begin{gathered}\gamma\in Q \\(\gamma|\bar{\Lambda}_{n-1})\geq 0, (\gamma|\bar{\Lambda}_1))\geq 0\end{gathered}}}(-1)^{l(w)}e^{wt_{q\gamma}(k\Lambda_0+\wr)-\wr}.	
\end{equation}
Furthermore, we have
\begin{enumerate}
	\item If $k=-3+\frac{2}{2m+1}$, $m\geq 1$, then  the maximal ideal of $V^k(sl_3)$ is generated by  two singular vectors $v^1$ and $v^2$ of weights $k\Lambda_0-3(2m+1)\delta+2\alpha_1+\alpha_2$ and $k\Lambda_0-3(2m+1)\delta+\alpha_1+2\alpha_2$, respectively. 
	
	\item If $k=-n+\frac{n-1}{q}$ satisfying $n\geq 4$, $q\geq 2$, and $(q,n-1)=1$, then the maximal ideal of $V^k(sl_n)$ is generated by a singular vector $w^1$ of weight $k\Lambda_0-2q\delta+\al_1+2\al_2+\cdots+2\al_{n-2}+\al_{n-1}$. 
\end{enumerate}
\end{theorem}
\begin{proof} Let $\g=sl_n$, $n\geq 3$.  It is obvious that $-\Lambda_0,  k\Lambda_0\in{\mathcal C}^+$,  for $k=-n+\frac{n-1}{q}$, $q\geq 2$, $(q,n-1)=1$. 
	
\vskip 0.2cm
Notice that if $q=1$, then $k=-1$ and $(\wW(-\Lambda_0), S(-\Lambda_0))=(\wW, S)$. So  by Lemma \ref{lweyl1}, $(\wW, S)\cong (\wW(k\Lambda_0), S(k\Lambda_0))$ through $\sigma$, for $k=-n+\frac{n-1}{q}$, $q\geq 2, (q,n-1)=1$.  
 In particular, 
\begin{equation}\label{emaximal24}
	\sigma (s_{\al_0})=s_{\beta_0}, \  \sigma(s_{\al_i})=s_{\al_i}, \ i=1,2,\cdots,n-1.	
\end{equation}	
Then for $w, y\in \wW$, 
\begin{equation}\label{emaximal9}
	w\geq_{-\Lambda_0}y \ {\rm if \ and \ if} \ \sigma(w)\geq_{k\Lambda_0}\sigma(y), 
\end{equation}
and 
\begin{equation}\label{emaximal10}
	l_{-\Lambda_0}(y)=l_{k\Lambda_0}\sigma(y), \ l_{-\Lambda_0}(w)=l_{k\Lambda_0}\sigma(w).
\end{equation}
\vskip 0.2cm
Let  $w\in\wW(-\Lambda_0)$ be such that $w$  is the longest element of $w\wW_0(-\Lambda_0)$. Then by  (\ref{emaximal24}) and (\ref{emaximal10}),  $\sigma(w)$ is the longest element of $\sigma(w)\wW_0(k\Lambda_0)$.  By Theorem \ref{KT}, 
\begin{equation}\label{emaximal11}
	{\rm ch}(L(w\circ (-\Lambda_0)))
	=\sum\limits_{w\leq_{-\Lambda_0}y\in \wW(-\Lambda_0)}(-1)^{l_{-\Lambda_0}(y)-l_{-\Lambda_0}(w)}Q^{-\Lambda_0}_{w,y}(1){\rm ch}(M(y\circ (-\Lambda_0))), 
\end{equation}
and	
\begin{equation}\label{emaximal25}
{\rm ch}(L(\sigma(w)\circ k\Lambda_0))
=\sum\limits_{\sigma(w)\leq_{k\Lambda_0}\sigma(y)\in \wW(k\Lambda_0)}(-1)^{l_{k\Lambda_0}(\sigma(y))-l_{k\Lambda_0}(\sigma(w))}Q^{k\Lambda_0}_{\sigma(w),\sigma(y)}(1){\rm ch}(M(\sigma(y)\circ k\Lambda_0)).
\end{equation}
Since $(\wW(-\Lambda_0), S(-\Lambda_0))\cong (\wW(k\Lambda_0), S(k\Lambda_0))$ via  the  automorphism $\sigma$, 
by the definition of $P^{\lambda}_{y,w}(q)$ and $Q^{\lambda}_{y,w}(q)$ \cite{KL79}, for $y,w\in\wW(-\Lambda_0)$,  
\begin{equation}\label{emaximal8}
P^{-\Lambda_0}_{y,w}(1)=P^{k\Lambda_0}_{\sigma(y),\sigma(w)}(1), \ Q^{-\Lambda_0}_{y,w}(1)=Q^{k\Lambda_0}_{\sigma(y),\sigma(w)}(1).
\end{equation}
Then together with (\ref{emaximal9})-(\ref{emaximal10}), we have
\begin{equation}\label{emaximal1}
	{\rm ch}(L(\sigma(w)\circ k\Lambda_0))=\sum\limits_{w\leq_{-\Lambda_0}y\in \wW(-\Lambda_0)}(-1)^{l_{-\Lambda_0}(y)-l_{-\Lambda_0}(w)}Q^{-\Lambda_0}_{w,y}(1){\rm ch}(M(\sigma(y)\circ k\Lambda_0)).
\end{equation}
The equation systems (\ref{emaximal11}) and  (\ref{emaximal1}) also implies that for $w,y\in\wW$,
\begin{equation}\label{emaximal22}
	[M(w\circ (-\Lambda_0)):L(y\circ (-\Lambda_0))]=[M(\sigma(w)\circ k\Lambda_0):L(\sigma(y)\circ k\Lambda_0)],
\end{equation}
where $[M: L(\mu)]$ denotes the multiplicity of the irreducible subquotient  $L(\mu)$ in $M$, for a module  $M$  in category ${\cal O}$.

\vskip 0.2cm
  By the character formulas given in \cite{KK79} for   generalized Verma modules of affine Lie algebras induced from integrable modules of the  associated simple Lie algebras, one can easily deduce   that  
\begin{equation}\label{max1}
[V(w\circ (-\Lambda_0)):L(y\circ (-\Lambda_0))]=[V(\sigma(w)\circ k\Lambda_0):L(\sigma(y)\circ k\Lambda_0)],
\end{equation}
where  $V(w\circ (-\Lambda_0))$ is the generalized  Verma module of the affine Lie algebra $\widehat{sl_n}$ induced from the finite-dimensional $sl_n$-module  $L(\overline{w\circ (-\Lambda_0}))$, and 
$V(\sigma(w)\circ (k\Lambda_0))$    is  the generalized  Verma module  induced from  $L(\overline{\sigma(w)\circ k\Lambda_0})$. 

\vskip 0.2cm
Notice that (\ref{emaximal1}) can also be described as
\begin{equation}\label{emaximal2}
	\widehat{R}{\rm ch}(L(\sigma(w)\circ k\Lambda_0))=\sum\limits_{w\leq_{-\Lambda_0}y\in \wW(-\Lambda_0)}(-1)^{l_{-\Lambda_0}(y)-l_{-\Lambda_0}(w)}Q^{-\Lambda_0}_{w,y}(1)e^{\sigma(y)\circ k\Lambda_0},
\end{equation}
where 
$$
\widehat{R}=\prod_{\al\in\wD_+}(1-e^{-\al})^{{\rm dim}{\fg}_{\al}}, \ \frak{g}=sl_n.
$$
By (\ref{emaximal26}), 
$$
\widehat{R}{\rm ch}L_{-1}(sl_n))=
\sum\limits_{w\in W}\sum\limits_{\tiny{\begin{gathered}\gamma\in Q \\(\gamma|\bar{\Lambda}_{n-1})\geq 0, (\gamma|\bar{\Lambda}_1))\geq 0\end{gathered}}}(-1)^{l(w)}e^{wt_{\gamma}(-\Lambda_0+\wr)-\wr}.
$$
By Lemma \ref{lweyl1} and Lemma \ref{lweyl2},  for $w\in W$, $\gamma\in Q$, 
$$
\sigma(wt_{\gamma})=wt_{q\gamma}.$$
Then we have 
$$
\widehat{R}{\rm ch}L_k(sl_n)=
\sum\limits_{w\in W}\sum\limits_{\tiny{\begin{gathered}\gamma\in Q \\(\gamma|\bar{\Lambda}_{n-1})\geq 0, (\gamma|\bar{\Lambda}_1))\geq 0\end{gathered}}}(-1)^{l(w)}e^{wt_{q\gamma}(k\Lambda_0+\wr)-\wr}.	
$$
(\ref{emaximal3}) is proved.

\vskip 0.2cm
\noindent
Set
$$
\Lambda^{(-1)}=\wW(-\Lambda_0)\circ(-\Lambda_0)=\{w\circ (-\Lambda_0)|w\in\wW(-\Lambda_0)\}$$
and
$$
\Lambda^{(k)}=\wW(k\Lambda_0)\circ k\Lambda_0=\{w\circ (k\Lambda_0)|w\in\wW(k\Lambda_0)\}.
$$
Let  ${\cal O}_{\Lambda^{(-1)}}$ and  $ {\cal O}_{\Lambda^{(k)}}$ be the corresponding indecomposable blocks,  respectively. 
Then by Theorem \ref{Fi06} and Lemma \ref{lweyl1},  

\begin{equation}\label{emaximal27}
	{\cal O}_{\Lambda^{(-1)}}\cong {\cal O}_{\Lambda^{(k)}}.
	\end{equation}
	We need to explain  the equivalence (\ref{emaximal27}) briefly.  Notice that  both $-\Lambda_0$ and $k\Lambda_0$ are  dominant. Then  for $w\in\wW(-\Lambda_0)$, $L(w\circ (-\Lambda_0))$ (resp. $L(\sigma(w)\circ k\Lambda_0)$) has a projective cover $P(w\circ (-\Lambda_0))$ (resp. $P(\sigma(w)\circ k\Lambda)$) in 	${\cal O}_{\Lambda^{(-1)}} $ (resp.  $ {\cal O}_{\Lambda^{(k)}}$).  Set
	  $$
	  {\cal P}_{\Lambda^{(-1)}}=\{P(w\circ (-\Lambda_0))|w\in\wW(-\Lambda_0)\}
	  $$
and
$$
{\cal P}_{\Lambda^{(k)}}=\{P(\sigma(w)\circ k\Lambda_0)|w\in\wW(-\Lambda_0)\}.
$$	  
$
{\cal P}_{\Lambda^{(-1)}}$ (resp.
${\cal P}_{\Lambda^{(k)}}$) is a faithful set of small projectives in the sense of \cite{Mit72}. View 
${\cal P}_{\Lambda^{(-1)}}$ (resp.
${\cal P}_{\Lambda^{(k)}}$)  as a full subcategory of ${\cal O}_{\Lambda^{(-1)}}$(resp. $ {\cal O}_{\Lambda^{(k)}}$). Let ${\mathbb C}-{\rm mod}^{{\cal P}_{\Lambda^{(-1)}}^{\rm opp}}$ (resp. ${\mathbb C}-{\rm mod}^{{\cal P}_{\Lambda^{(k)}}^{\rm opp}}$) be the category of additive functors ${\cal P}_{\Lambda^{(-1)}}\rightarrow {\C}-{\rm mod}$ (resp. ${\cal P}_{\Lambda^{(k)}}\rightarrow {\C}-{\rm mod}$). By Theorem 3.1 of \cite{Mit72},
$$
	{\cal O}_{\Lambda^{(-1)}}\rightarrow {\mathbb C}-{\rm mod}^{{\cal P}_{\Lambda^{(-1)}}^{\rm opp}} \ \ ({\rm resp.}  \ 
		{\cal O}_{\Lambda^{(k)}}\rightarrow {\mathbb C}-{\rm mod}^{{\cal P}_{\Lambda^{(k)}}^{\rm opp}} )
$$
$$
	M\rightarrow {\rm Hom}(\cdot,  \ M) \ \  \ ({\rm resp.} \	\tilde{M}\rightarrow  {\rm Hom}(\cdot, \  \tilde{M})).
$$
is an equivalence of ${\C}$-categories.  Theorem \ref{Fi06} (i.e., Theorem 11 of   \cite{Fi06})  tells us  that the Coxeter group automorphism $\sigma$ determines an  equivalence between  the categories ${\cal O}_{\Lambda^{(-1)}}$ and ${\cal O}_{\Lambda^{(k)}}$  via  the above $\C$-category equivalences (c.f.  \cite{Fi06} for details). 
\vskip 0.2cm	  
	    Let $F$ be the covariant  functor   from 
	  $	{\cal O}_{\Lambda^{(-1)}}$ to ${\cal O}_{\Lambda^{(k)}}$. 
	   Then in particular, $F$  maps the class of Verma module $M(w\circ (-\Lambda_0))$ to  the class of Verma module $M(\sigma(w)\circ k\Lambda_0)$, the class of simple module $L(w\circ (-\Lambda_0))$ to the class of simple module $L(\sigma(w)\circ k\Lambda_0)$, and  it is easy to see that $F$ maps the class of  $V(w\circ(-\Lambda_0))$
	  to  the class of $V(\sigma(w)\circ k\Lambda_0)$, 
	   for $w\in\wW(-\Lambda_0)$.

\vskip 0.2cm
{\bf Case 1}  \  $n=3$,  $k=-3+\frac{2}{2m+1}$, $m\geq 1$.  Let $u^1$ and $u^2$ be the singular vectors given in (1) of Theorem \ref{AP08} with highest weights $-\Lambda_0-3\delta+2\al_1+\al_2$ and $-\Lambda_0-3\delta+\al_1+2\al_2$, respectively. It is easy to check that 
\begin{equation}\label{emaximal4}
s_{\al_2}\circ t_{\al_1}\circ(-\Lambda_0)=s_{\al_2}\circ (t_{\al_1}(-\Lambda_0+\wr)-\wr)=-\Lambda_0-3\delta+2\al_1+\al_2,  
\end{equation}
and
\begin{equation}\label{emaximal5}
s_{\al_1}\circ t_{\al_2}\circ(-\Lambda_0)=s_{\al_1}\circ (t_{\al_2}(-\Lambda_0+\wr)-\wr)=-\Lambda_0-3\delta+\al_1+2\al_2.
\end{equation}
Then by (\ref{emaximal22}) and Lemma \ref{hom}, 
$V^k(sl_3)$ has two singular vectors $v^1$ and $v^2$ of weights $s_{\al_2}\circ t_{(2m+1)\al_1}\circ k\Lambda_0$ and $s_{\al_1}\circ t_{(2m+1)\al_2}\circ k\Lambda_0$, respectively.  It is easy to check  that
\begin{equation}\label{emaximal6}
s_{\al_2}\circ t_{(2m+1)\al_1}\circ k\Lambda_0=s_{\al_2}\circ (t_{(2m+1)\al_1}(k\Lambda_0+\wr)-\wr)=k\Lambda_0-3(2m+1)\delta+2\al_1+\al_2,
\end{equation}and
\begin{equation}\label{emaximal7}
s_{\al_1}\circ t_{(2m+1)\al_2}\circ k\Lambda_0=s_{\al_1}\circ (t_{(2m+1)\al_2} (k\Lambda_0+\wr)-\wr)=k\Lambda_0-3(2m+1)\delta+\al_1+2\al_2.
\end{equation}
Since 
$$
[V^{-1}(sl_3): L(s_{\al_1}t_{\al_2}\circ(-\Lambda_0))]=[V^{-1}(sl_3): L(s_{\al_2}t_{\al_1}\circ(-\Lambda_0))]=1,
$$
by  (\ref{max1}), 
\begin{equation}\label{emax6}
[V^{k}(sl_3): L(s_{\al_1}t_{(2m+1)\al_2}\circ k\Lambda_0)]=[V^{k}(sl_3): L(s_{\al_2}t_{(2m+1)\al_1}\circ k\Lambda_0)]=1.
\end{equation}
Let $\tilde{I}$ be the ideal of $V^k(sl_3)$ generated by $v^1$ and $v^2$.  
Since $I$ is the maximal ideal of $V^{-1}(sl_3)$,  we have the following exact sequence:
$$
0\rightarrow I \rightarrow V^{-1}(sl_3)\rightarrow L_{-1}(sl_3)\rightarrow 0.
$$
Then we have the exact sequence:
$$
0\rightarrow F(I) \rightarrow V^{k}(sl_3)\rightarrow L_{k}(sl_3)\rightarrow 0,
$$
since $ F(V^{-1}(sl_3))=V^{k}(sl_3)$, $F(L_{-1}(sl_3))=L_{k}(sl_3)$. 

\vskip 0.2cm
On the other hand, since   $I$ is a quotient of $M(s_{\al_1}\circ t_{\al_2}\circ(-\Lambda_0))\oplus M(s_{\al_2}\circ t_{\al_1}\circ(-\Lambda_0))$, it  follows that   $F(I)$ is a quoitent of $M(s_{\al_1}\circ t_{(2m+1)\al_2}\circ k\Lambda_0)\oplus M(s_{\al_2}\circ t_{(2m+1)\al_1}\circ k\Lambda_0)$. This together with (\ref{emax6}) means that $F(I)$ 
 must be 
 the ideal $\tilde{I}$  generated by $v^1$ and $v^2$. 
 Since $F(I)$ is maximal, it follows  that $\tilde{I}$ is maximal.

\vskip 0.2cm
{\bf Case 2} \ $n\geq 4$, $k=-n+\frac{n-1}{q}$, $q\geq 2$, $(q, n-1)=1$.  Let $u$ be the singular vector in (2) of Theorem \ref{AP08}. The weight of $u$ with respect to the Cartan subalgebra $\fh$ is 
$-\Lambda_0-2\delta+\theta+\beta$, where $\beta=\al_2+\cdots +\al_{n-2}$. It is easy to  check that
$$
-\Lambda_0-2\delta+\theta+\beta=(s_{\theta}s_{\al_1}s_{\al_{n-1}})\circ (t_{-\beta}(-\Lambda_0+\wr)-\wr).
$$
Then by (\ref{emaximal22}) and Lemma \ref{hom}, $V^k(sl_n)$ has a singular vector $w^1$ of weight $(s_{\theta}s_{\al_1}s_{\al_{n-1}}t_{-q\beta})\circ k\Lambda_0$.
One can easily check that 
$$
(s_{\theta}s_{\al_1}s_{\al_{n-1}})\circ (t_{-q\beta}(k\Lambda_0+\wr)-\wr)=k\Lambda_0-2q\delta+\theta+\beta.
$$
By (4) of Theorem \ref{AP08}, the ideal $J$  generated by $u$ is the maximal ideal of $V^{-1}(sl_n)$.   Then we have the exact sequence:
\begin{equation}\label{emax7}
0\rightarrow J\rightarrow V^{-1}(sl_n)\rightarrow L_{-1}(sl_n)\rightarrow 0.
\end{equation}
By the fact that   
$F(V^{-1}(sl_n))=V^k(sl_n)$, and $F(L_{-1}(sl_n))=L_k(sl_n)$,
we have 
\begin{equation}\label{emax8}
0\rightarrow F(J) \rightarrow V^{k}(sl_n)\rightarrow L_{k}(sl_n)\rightarrow 0.
\end{equation}
Notice that $J$ is a quotient of $M(s_{\theta}s_{\al_1}s_{\al_{n-1}}t_{-\beta})\circ (-\Lambda_0)) $,  it follows that $F(J)$ is a quotient of  $M((s_{\theta}s_{\al_1}s_{\al_{n-1}}t_{-q\beta})\circ k\Lambda_0)$. Then  $F(J)$ is  generated by $w^1$, since 
$$[V^{-1}(sl_n): L((s_{\theta}s_{\al_1}s_{\al_{n-1}}t_{-\beta})\circ (-\Lambda_0))]=[V^{k}(sl_n): L((s_{\theta}s_{\al_1}s_{\al_{n-1}}t_{-q\beta})\circ k\Lambda_0)]=1.$$
 By (\ref{emax8}),  $F(J)$  is maximal. We complete the proof of (2).  We remark   that  (\ref{emax8}) can also be deduced  from (\ref{emax7}) and (\ref{max1}) without using the category equivalence (\ref{emaximal27}).
\end{proof}

\section{Varieties of $L_k(sl_3)$ for $k=-3+\frac{2}{2m+1}$, $m\geq 0$}
Let $\g$ and $V^k(\g)$ be as above in Section 3. 
For $x\in V^k(\frak{g})$, we denote by $\bar{x}$ the image of $x$ in $R_{V^k(\g)}$. Recall that   $R_{V^k(\g)}\cong \C[\g^*]$ by
$$
\overline {x_1(-1)\cdots x_n(-1){\bf 1}}\mapsto x_1\cdots x_n,
$$
for $x_1,\cdots,x_n\in \g$. 
As in Section 2, we denote by ${\mathcal I}_k$  the maximal ideal of $V^k(\g)$, then 
$$
L_k(\g)=V^k(\g)/{\mathcal I}_k
$$
and
$$
R_{L_k(\g)}=\C[\g^*]/\overline{{\mathcal I}_k},
$$
where $\overline{{\mathcal I}_k}$ is the image of $\mathcal{I}_k$ in $R_{V^k(\g)}$. 
\subsection{The associated variety of $L_{-1}(sl_3)$}
Let $\Delta$ be the root system of $sl_3({\mathbb C})$ and 
$$\{h_1, h_2, e_{\al_1}, e_{\al_2}, e_{\al_1+\al_2}, f_{\al_1},f_{\al_2}, f_{\al_1+\al_2}\}
$$
be a Chevalley basis of $sl_3({\mathbb C})$ with structure constants $c_{\al,\be}$, for $\al,\be\in\Delta$ such that $$c_{\al_1,\al_2}=1.$$
Then we have 
$$
c_{\al_2,-\al_1-\al_2}=c_{-\al_1-\al_2,\al_1}=1.
$$
\vskip 0.2cm
\noindent
Recall from Theorem \ref{AP08} \cite{AP08} that 
$$
\begin{array}{ll}
u^1=& [-e_{\al_1}(-1)^2e_{\al_2}(-1)+e_{\al_1}(-1)e_{\al_1+\al_2}(-1)h_2(-1)+e_{\al_1+\al_2}(-1)^2f_{\al_2}(-1)]{\bf 1},\\\\
u^2=&[e_{\al_1}(-1)e_{\al_2}(-1)^2+e_{\al_2}(-1)e_{\al_1+\al_2}(-1)h_1(-1)\\
&-2e_{\alpha_2}(-1)e_{\al_1+\al_2}(-2)
-e_{\al_1+\al_2}(-1)^2f_{\al_1}(-1)]{\bf 1}.
\end{array}
$$
are  singular vectors of $V^{-1}(sl_3)$. Then  we have the following lemma.
\begin{lem}\label{l3.2} There exist $x,y\in {\mathcal I}_{-1}$ such that
$$
\begin{array}{ll}
x\equiv & -e_{\al_1+\al_2}(-1)h_1(-1)h_2(-1){\bf 1}-e_{\al_1}(-1)e_{\al_1+\al_2}(-1)f_{\al_1}(-1){\bf 1}\\\\
&-e_{\al_1+\al_2}(-1)^2f_{\al_1+\al_2}(-1){\bf 1}
+e_{\al_1}(-1)e_{\al_2}(-1)(2h_1+h_2)(-1){\bf 1}\\\\
&+2e_{\al_2}(-1)e_{\al_1+\al_2}(-1)f_{\al_2}(-1){\bf 1}
(\mod~ F^1(V^{-1}(sl_3)))
\\\\
y\equiv & h_1(-1)h_2(-1)(h_1(-1)+h_2(-1)){\bf 1}+e_{\al_1}(-1)e_{\al_2}(-1)f_{\al_1+\al_2}(-1){\bf 1}\\\\
&-3e_{\al_1+\al_2}(-1)f_{\al_1}(-1)f_{\al_2}(-1){\bf 1}-e_{\al_1}(-1)f_{\al_1}(-1)h_1(-1){\bf 1}\\\\
& +e_{\al_1+\al_2}(-1)f_{\al_1+\al_2}(-1)(h_1(-1)+h_2(-1)){\bf 1} \\\\ &-e_{\al_2}(-1)f_{\al_2}(-1)h_2(-1){\bf 1}\ (\mod~ F^1(V^{-1}(sl_3))).
\end{array}
$$
\end{lem}
\begin{proof} Let $x=f_{\al_1}(0)u^1$, and   $y=f_{\al_1}(0)f_{\al_1+\al_2}(0)u^1\in {\mathcal I}_{-1}$. Then the lemma follows.	
	\end{proof}
Let 
$$
\begin{array}{ll}
p_1=& -e_{\al_1}^2e_{\al_2}+e_{\al_1}e_{\al_1+\al_2}h_2+e_{\al_1+\al_2}^2f_{\al_2},
\\\\
p_2=& e_{\al_1}e_{\al_2}^2+e_{\al_2}e_{\al_1+\al_2}h_1-e_{\al_1+\al_2}^2f_{\al_1},\\\\
p_3= &-e_{\al_1+\al_2}h_1h_2-e_{\al_1}e_{\al_1+\al_2}f_{\al_1}-e_{\al_1+\al_2}^2f_{\al_1+\al_2}
 +e_{\al_1}e_{\al_2}(2h_1+h_2)
+2e_{\al_2}e_{\al_1+\al_2}f_{\al_2},
\\\\
p_4=&h_1h_2(h_1+h_2)+e_{\al_1}e_{\al_2}f_{\al_1+\al_3}
-3e_{\al_1+\al_2}f_{\al_1}f_{\al_2}-e_{\al_1}f_{\al_1}h_1,\\\\
&
+e_{\al_1+\al_2}f_{\al_1+\al_2}(h_1+h_2)-e_{\al_2}f_{\al_2}h_2.
\end{array}
$$
Then by Lemma \ref{l3.2},  $p_i\in \overline{{\mathcal I}_{-1}}$, $1\leq i\leq 4$. We have the following result which is also covered in \cite{AFK24} by different method.
\noindent
\begin{theorem}\label{variety1}
	 Denote $\lambda=h_1-h_2$. Then 
$$
X_{L_{-1}(sl_3)}=\overline{{\mathbb S}_{\frak l_1}}= \overline{G.{\mathbb C}^*\lambda}=G.{\mathbb C}^*\lambda\cup\overline{{\mathcal O}_{\min}},
$$
which is the closure of the Dixmier sheet ${\mathbb S}_{\frak l_1}=G.{\mathbb C}^*\lambda$, where ${\frak l}_1=\frak{h}\oplus \C e_{\theta}\oplus \C f_{\theta}$.  In particular, $X_{L_{-1}(sl_3)}$ is irreducible and 
$$
\dim X_{L_{-1}(sl_3)}=5,
$$
where ${\mathcal O}_{\min}$ is the minimal non-zero nilpotent orbit of $sl_3({\mathbb C})$.
\end{theorem}
\begin{proof}  It is known that for any simple Lie algebra $\g$ over ${\C}$, $\g$ is the finite disjoint union of its Jordan classes \cite{TY05}.  In particular, for $\g=sl_3(\C)$, we have 
	$$\g=G\cdot (\C^*\lambda)\cup G\cdot (\C^*\lambda+f_{\theta})\cup G\cdot ({\frak h}^{reg})\cup G\cdot f_{\theta}\cup G\cdot (f_{\al_1}+f_{\al_2})\cup\{0\},$$
	where $G$ is the connected adjoint group of $\g$, and $\theta=\al_1+\al_2$.	
	Notice that
	 $$p_1(f_{\al_1}+f_{\al_2})=-1, \ p_3(t\lambda+f_{\theta})=-t^2.$$  
	 This implies that  $G.(f_{\al_1}+f_{\al_2})$ and $G.(\C^*\lambda+f_{\theta})$ could not belong in $X_{L_{-1}(sl_3)}$. 
Furthermore, for regular semisimple vector $h$, all
$
(h_1|h), \ (h|h_2), (h|h_1+h_2)
$ are  non-zero. 
So
$$
p_4(h)=(h|h_1)(h|h_2)(h|h_1+h_2)\neq 0.
$$ We see that for any $h\in {\frak h}^{reg}$, $h\notin X_{L_{-1}(sl_3)}$. 
We deduce that
$$
X_{L_{-1}(sl_3)}\subseteq \overline{G.{\mathbb C}^*\lambda}=G.{\mathbb C}^*\lambda\cup\overline{{\mathcal O}_{\min}}.
$$
By (3) of Theorem \ref{AP08},  the ideal generated by $u^1$ and $u^2$ is the maximal ideal of $V^{-1}(sl_3)$. This means that 
$$
X_{L_{-1}(sl_3)}= \overline{G.{\mathbb C}^*\lambda}=G.{\mathbb C}^*\lambda\cup\overline{{\mathcal O}_{\min}}.
$$
It follows that $X_{L_{-1}(sl_3)}$ is irreducible and  $
\dim X_{L_{-1}(sl_3)}=5$. 
\end{proof}
\subsection{Associated varieties of $L_{k}(sl_3)$ for $k=-3+\frac{2}{2m+1}, m\geq 1$}
In this subsection we always  assume that $k=-3+\frac{2}{2m+1}$, $m\geq 1$.

\vskip 0.2cm
Recall from Theorem \ref{tmain3.1} that 
		$V^k(sl_3)$ has two singular vectors $v^1$ and $v^2$ with weights $k\Lambda_0-3(2m+1)\delta+2\alpha_1+\alpha_2$ and $k\Lambda_0-3(2m+1)\delta+\alpha_1+2\alpha_2$, respectively.  
  Then for $i=1,2$, 
$$
e_{\al_1}(0)v^i=e_{\al_2}(0)v^i=f_{\al_1}(1)v^i=f_{\al_2}(1)v^i=0.
$$
Notice that each element in $V^k(sl_3)$ is a linear combination of elements of the following form:
 \begin{align}
\label{eq:PBW_basis}
z = z^{(+)}  z^{(-)}  z^{(0)}  {\bf 1},
\end{align}
with
\begin{align*}
 z^{(+)} := & e_{\al_{1}}(-1)^{a_{11}}\cdots
e_{\al_{1}} (- r_1)^{a_{1r_1}}e_{\al_{2}}(-1)^{a_{21}}\cdots
e_{\al_{2}} (- r_2)^{a_{2r_2}}\\ \nonumber
& e_{\al_1+\al_2}(- 1)^{a_{31}}\cdots e_{\al_1+\al_2}(- r_3)^{a_{3r_3}} , \\ \nonumber
 z^{(-)} : =& f_{\al_1}(-1)^{b_{11}}\cdots f_{\al_1}(- s_1)^{b_{1s_1}}f_{\al_2}(-1)^{b_{21}}\cdots f_{\al_2}(- s_2)^{b_{2s_2}}\\ \nonumber
& f_{\al_1+\al_2}(- 1)^{b_{31}}
\cdots f_{\al_1+\al_2}(- s_3)^{b_{3s_3}} , \\\nonumber
  z^{(0)} := & h_1(- 1)^{c_{11}}\cdots h_1(- t_1)^{c_{1t_1}}
h_2(- 1)^{c_{21}}\cdots h_2(- t_2)^{c_{2t_2}},
\end{align*}
where $r_1,r_2,r_3, s_1,s_2,s_3,,t_1,t_2$
are positive integers, and $a_{lp},b_{ln},c_{ij}$, for $l =1,2,3$,
$p=1,\ldots,r_l$, $n=1,\ldots,s_l$,  $i = 1, 2$, $j=1,\ldots,t_i$, 
are nonnegative integers such
that at least one of them is non-zero. Recall  that
$$
depth(z^{(+)})=\sum\limits_{i=1}^3\sum\limits_{j=1}^{r_{i}-1}(j-1)a_{ij}, ~ depth(z^{(-)})=\sum\limits_{i=1}^3\sum\limits_{j=1}^{s_{i}-1}(j-1)b_{ij}, $$
$$ depth(z^{(0)})=\sum\limits_{i=1}^2\sum\limits_{j=1}^{t_{i}-1}(j-1)c_{ij};
$$
$$
deg{(z^{(+)})}=\sum\limits_{i=1}^3\sum\limits_{j=1}^{r_{i}}a_{ij}, ~ deg{(z^{(-)})}=\sum\limits_{i=1}^3\sum\limits_{j=1}^{s_{i}}b_{ij},~ deg{(z^{(0)})}=\sum\limits_{i=1}^2\sum\limits_{j=1}^{t_{i}}c_{ij};
$$
$$
depth(z)=depth(z^{(+)})+depth(z^{(-)})+depth(z^{(0)});$$
$$ deg(z)=deg(z^{(+)})+deg(z^{(-)})+deg(z^{(0)}).
$$

 Let $V^1$ be the subspace of $V^k(sl_3)$ linearly spanned by elements $z=z^{(+)}  z^{(-)}  z^{(0)}  {\bf 1}$ of weight $k\Lambda_0-3(2m+1)\delta+2\al_1+\al_2$ such that $deg(z^{(0)})\leq 6m-2$ or $depth(z^{(0)})\geq 1$.  Then we  may assume that 
$$
\begin{array}{ll}
v^1=&\sum\limits_{i=0}^{6m+1}a_ie_{\al_1}(-1)e_{\al_1+\al_2}(-1)h_1(-1)^ih_2(-1)^{6m+1-i}{\bf 1}\\\\
&+\sum\limits_{i=0}^{6m}[x_ie_{\al_1}(-1)^2e_{\al_2}(-1)+y_ie_{\al_1+\al_2}(-1)^2f_{\al_2}(-1)+b_ie_{\al_1}(-2)e_{\al_1+\al_2}(-1)\\\\
& ~~~+c_ie_{\al_1}(-1)e_{\al_1+\al_2}(-2)]h_1(-1)^ih_2(-1)^{6m-i}{\bf 1}\\\\
&+\sum\limits_{i=0}^{6m-1}[d_ie_{\al_1+\al_2}(-1)e_{\al_1+\al_2}(-2)f_{\al_2}(-1)
+z_ie_{\al_1}(-1)e_{\al_2}(-1)e_{\al_1+\al_2}(-1)\\\\
&f_{\al_2}(-1)
+l_ie_{\al_1}(-1)^2e_{\al_1+\al_2}(-1)f_{\al_1}(-1)
+n_ie_{\al_1}(-2)e_{\al_1+\al_2}(-2)\\\\
&+k_ie_{\al_1}(-1)^2e_{\al_2}(-2)+g_ie_{\al_1}(-1)e_{\al_1}(-2)e_{\al_2}(-1)+p_ie_{\al_1+\al_2}(-1)^2f_{\al_2}(-2)\\\\
&+q_ie_{\al_1}(-1)e_{\al_1+\al_2}(-1)^2f_{\al_1+\al_2}(-1)+m_ie_{\al_1}(-1)e_{\al_1+\al_2}(-3)\\\\
&+r_ie_{\al_1}(-3)e_{\al_1+\al_2}(-1)]h_1(-1)^ih_2(-1)^{6m-1-i}{\bf 1}+u^1,
\end{array}
$$
where $u^1\in V^1$. 

\vskip 0.3cm
The  following result comes from Lemma 3.1 of  \cite{AJM21}. 
\begin{lem}\label{AJM21} Not all $a_i, x_j, 0\leq i\leq 6m+1, 0\leq j\leq 6m$ in $v^1$ are zero. 
	
	\end{lem}

We now have the following lemma.
\begin{lem}\label{l3.2.2} 
\begin{equation}\label{e3.2.28}
a_i-y_i=0,  \ a_{6m+1}=0,  \ 0\leq i\leq 6m,
\end{equation}
\begin{equation}\label{e3.2.3}
-2(i+1)a_{i+1}+(6m+1-i)a_i+x_i+l_{i-1}=0,  \ 0\leq i\leq 6m,
\end{equation}
\begin{equation}\label{e3.2.11}
(i+1)a_{i+1}-2(6m+1-i)a_i-2x_i+z_i=0, \ 0\leq i\leq 6m-1,
\end{equation}
\begin{equation}\label{e3.2.19}
(6m+1)a_{6m+1}-2a_{6m}-2x_{6m}=0,
\end{equation}
\begin{equation}\label{e3.2.6}
(-4+\frac{2}{2m+1})a_i+2y_i-b_{i-1}=0, \ \ 0\leq i\leq 6m,
\end{equation}
\begin{equation}\label{e3.2.13}
(-4+\frac{2}{2m+1})a_{6m+1}-b_{6m}=0, \   (-4+\frac{2}{2m+1})a_0+2y_0=0,
\end{equation}
\begin{equation}\label{e3.2.12}
(-4+\frac{2}{2m+1})a_i-2x_i-c_{i-1}-c_i=0, \ 0\leq i\leq 6m,
\end{equation}
\begin{equation}\label{e3.2.14}
a_0=y_0=0,
\end{equation}
\begin{equation}\label{e3.2.4}
a_{i}+x_i=0, \ 0\leq i\leq 6m,
\end{equation}
\begin{equation}\label{e3.2.4c}
-(i+1)a_{i+1}-(6m-i)a_i+q_{i}+q_{i-1}=0, \ 0\leq i\leq 6m.
\end{equation}
	\end{lem}
\begin{proof}
 We consider $e_{\al_1}(0)v^1$, $e_{\al_2}(0)v^1$, $f_{\al_1}(1)v^1$, and  $f_{\al_1+\al_2}(1)v^1$. By the definition of $V^1$, there are no monomials 
$z=z^{(+)}  z^{(-)}  z^{(0)}  {\bf 1}$   in  $e_{\al_1}(0)u^1$, $e_{\al_2}(0)u^1$, $f_{\al_1}(1)u^1$, $f_{\al_2}(1)u^1$,  and  $f_{\al_1+\al_2}(1)u^1$ such that $deg(z^{(0)})\geq 6m$, $depth(z)=0$, 
and $z^{(-)}=1$. Then it is easy to deduce that 
the coefficients of $e_{\al_1}(-1)^2e_{\al_1+\al_2}(-1)h_1(-1)^ih_2(-1)^{6m-i}{\bf 1}$ in $e_{\al_1}(0)v^1$, $0\leq i\leq 6m$,   are 
$$
-2(i+1)a_{i+1}+(6m+1-i)a_i+x_i+l_{i-1}
,  \ 0\leq i\leq 6m.
$$
Then (\ref{e3.2.3}) holds.

\vskip 0.2cm
\noindent
 The coefficients of $e_{\al_1+\al_2}(-1)^2h_1(-1)^ih_2(-1)^{6m+1-i}{\bf 1}$ in $e_{\al_2}(0)v^1$, $0\leq i\leq 6m$,  are 
$$
-a_i+y_i, \  0\leq i\leq 6m,
$$
and the coefficient of $e_{\al_1+\al_2}(-1)^2h_1(-1)^{6m+1-i}{\bf 1}$ is $a_{6m+1}$. Then  (\ref{e3.2.28}) follows. 

\vskip 0.2cm
\noindent
The coefficients of $e_{\al_1}(-1)e_{\al_2}(-1)e_{\al_1+\al_2}(-1)h_1(-1)^ih_2(-1)^{6m-i}{\bf 1}$, $0\leq i\leq 6m-1$, in $e_{\al_2}(0)v^1$ are
$$
(i+1)a_{i+1}-2(6m+1-i)a_i-2x_i+z_i, \ 0\leq i\leq 6m-1,
$$
and the coefficient  of $e_{\al_1}(-1)e_{\al_2}(-1)e_{\al_1+\al_2}(-1)h_1(-1)^{6m}{\bf 1}$  in $e_{\al_2}(0)v^1$ is
$$
(6m+1)a_{6m+1}-2a_{6m}-2x_{6m}.
$$
Then (\ref{e3.2.11}) and (\ref{e3.2.19}) hold. 

\vskip 0.2cm
\noindent
The coefficients of $e_{\al_1+\al_2}(-1)h_1(-1)^ih_2(-1)^{6m+1-i}{\bf 1}$, $0\leq i\leq 6m+1$,
 in $f_{\al_1}(1)v^1$ are
 $$
  (-4+\frac{2}{2m+1})a_i+2y_i-b_{i-1}, \ \ 1\leq i\leq 6m,
 $$
 and
 $$
  (-4+\frac{2}{2m+1})a_{6m+1}-b_{6m}(i=6m+1), \   (-4+\frac{2}{2m+1})a_0+2y_0(i=0).
 $$
 Then (\ref{e3.2.6})and (\ref{e3.2.13}) hold.

\vskip 0.2cm
\noindent
The coefficients of $e_{\al_1}(-1)h_1(-1)^ih_2(-1)^{6m+1-i}{\bf 1}$, $0\leq i\leq 6m$, in $f_{\al_1+\al_2}(1)v^1$ are 
$$
(-4+\frac{2}{2m+1})a_i-2x_i-c_{i-1}-c_i, \ 0\leq i\leq 6m.
$$
This implies (\ref{e3.2.12}).

\vskip 0.2cm
\noindent
By (\ref{e3.2.4}), (\ref{e3.2.13}), and the fact that $m\geq 1$, 
$$
a_0=y_0=0.
$$
Then (\ref{e3.2.14}) follows.

\vskip 0.2cm
\noindent
Notice that $(\al_2|2\al_1+\al_2)=0$ and $e_{\al_2}(0)v^1=0$. It follows that 
$f_{\al_2}(0)v^1=0$. 
Considering the coefficients of $e_{\al_1}(-1)^2h_1(-1)^ih_2(-1)^{6m+1-i}$, $0\leq i\leq 6m$, we obtain that 
$$
a_i+x_i=0, \ 0\leq i\leq 6m,
$$
which is  (\ref{e3.2.4}). 

\vskip 0.2cm
\noindent
The coefficients of $e_{\al_1}(-1)e_{\al_1+\al_2}(-1)^2h_1(-1)^ih_{2}(-1)^{6-i}{\bf 1}$, $0\leq i\leq 6m$, in $e_{\al+\al_2}(0)v^1$ are
$$
-(i+1)a_{i+1}-(6m+1-i)a_i+y_i+q_{i}+q_{i-1}, \ 0\leq i\leq 6m.
$$
Then (\ref{e3.2.4c}) follows from (\ref{e3.2.28}).
\end{proof}
By Lemma  \ref{AJM21}, not all $a_i, x_j, 0\leq i\leq 6m+1, 0\leq j\leq 6m$ are zero, and by Lemma \ref{l3.2.2}, $a_i=-x_i, 0\leq i\leq 6m, \ a_{6m+1}=0$. Then  not all $a_i, 0\leq i\leq 6m+1$ are zero. By (\ref{e3.2.4}), $a_0=0$. Thus we may assume that 
$$
a_0=\cdots=a_{r-1}=0, \ a_r\neq 0,
$$
for some $r\geq 1$. 
We have the following lemma.
\begin{lem}\label{l3.2.6} For $k=-3+\frac{2}{2m+1}$,  $m\geq 1$, 
$
	r=2m.
	$
	\end{lem}
\begin{proof}  
	It can be checked directly that the coefficient of $$e_{\al_1}(-1)e_{\al_1+\al_2}(-1)h_{1}(-1)^{r-1}h_2(-1)^{6m+1-r}$$ in $h_1(1)v^1$
	is 
\begin{equation}\label{e3.2.2c}	
-5ra_r+2r(-3+\frac{2}{2m+1})a_r+4x_{r-1}+2y_{r-1}+2b_{r-1}+c_{r-1}+c_{r-2}-z_{r-1}+4l_{r-2}+2(q_{r-1}+q_{r-2}),
\end{equation}	
which should be zero.
By (\ref{e3.2.28}) and (\ref{e3.2.4}), 
\begin{equation}\label{e3.2.1c}
x_{r-1}=y_{r-1}=0.
\end{equation}
Then by (\ref{e3.2.3}) and (\ref{e3.2.11}), we have 
$$
l_{r-2}=2ra_r, \ z_{r-1}=-ra_r.
$$
By (\ref{e3.2.28}) and (\ref{e3.2.6}),
$$
b_{r-1}=\frac{-4m}{2m+1}a_r.
$$
By (\ref{e3.2.12}) and (\ref{e3.2.1c}),
$$
c_0=0, \ c_j+c_{j-1}=0,  1\leq j\leq r-1.
$$
Then we have
$$
c_0=\cdots=c_{r-1}=0.
$$
By (\ref{e3.2.4c}), 
$$
q_{r-1}+q_{r-2}=ra_r.
$$
Thus (\ref{e3.2.2c}) becomes
$$
(\frac{4r}{2m+1}-\frac{8m}{2m+1})a_r=0.
$$
This deduces that $r=2m$.
\end{proof}
Let $U^0$ be the $sl_3$-module generated by $v^1$. The following lemma is easy to check. 
\begin{lem}\label{l3.2.7}  
	With respect to the Cartan subalgebra $\frak{h}=\C h_1\oplus\C h_2$,  the weight zero subspace of $U^0$ is one-dimensional, which is linearly spanned by $f_{\al_1}(0)f_{\al_1+\al_2}(0)v^1$. 
	\end{lem}
Let $W^1$ be the subspace of $V^k(sl_3)$ linearly spanned by monomials $z=z^{(+)}z^{(-)}z^{(0)}{\bf 1}$ such that $z\in F^0(V^k(sl_3))$ and $z^{(-)}\neq 1$.  Then  it can be  checked directly  that
\begin{equation}\label{e3.2.30}
f_{\al_1}(0)f_{\al_1+\al_2}(0)v^1
=\sum\limits_{i=2m}^{6m+1}a_ih_1(-1)^{i+1}h_2(-1)^{6m+1-i}(h_1+h_2)(-1){\bf 1}+w^1, 
\end{equation}
for some $w^1\in W^1+F^1(V^k(sl_3))$. We denote
$$
v^3=f_{\al_1}(0)f_{\al_1+\al_2}(0)v^1,
$$
and without loss of generality we may assume that 
$$
a_{2m}=1.
$$
We have the following lemma.
\begin{lem}\label{l3.2.8}
	\begin{equation}\label{e3.2.5c}
	\sum\limits_{i=2m}^{6m+1}a_ih_1(-1)^{i+1}h_2(-1)^{6m+1-i}(h_1+h_2)(-1){\bf 1}=h_1(-1)^{2m+1}h_2(-1)^{2m+1}(h_1+h_2)(-1)^{2m+1}{\bf 1}.
	\end{equation}	
That is,	
$$
\begin{array}{ll}
&v^1-e_{\al_1}(-1)e_{\al_1+\al_2}(-1)h_1(-1)^{2m}h_2^{2m+1}(h_1+h_2)(-1)^{2m}{\bf 1}\\\\
&+e_{\al_1}(-1)^2e_{\al_2}(-1)h_1(-1)^{2m}h_2(-1)^{2m}(h_1+h_2)(-1)^{2m}
{\bf 1}\in W^1+F^1(V^k(sl_3)).
\end{array}
$$	
and	
\begin{equation}\label{e.3.2.3c}
v^3=h_1(-1)^{2m+1}h_2(-1)^{2m+1}(h_1+h_2)(-1)^{2m+1}{\bf 1}+w^1,
\end{equation}
where $w^1$ is the same as in (\ref{e3.2.30}).
\end{lem}
\begin{proof}
	It suffices to prove that 
$$
\sum\limits_{i=2m}^{6m+1}a_ih_1(-1)^{i+1}h_2(-1)^{6m+1-i}(h_1+h_2)(-1){\bf 1}=h_1(-1)^{2m+1}h_2(-1)^{2m+1}(h_1+h_2)(-1)^{2m+1}{\bf 1}.
$$
Denote
$$
\pi_2=\exp f_{\al_2}(0)\exp(-e_{\al_2}(0))\exp f_{\al_2}(0), $$$$ \pi_3
=\exp f_{\al_1+\al_2}(0)\exp(-e_{\al_1+\al_2}(0))\exp f_{\al_1+\al_2}(0).$$ 
Then  by Lemma \ref{l3.2.7}, there exists $0\neq c_i\in \C$, $i=2,3$ such that 
$$
\pi_i(v^3)=c_iv^3.
$$
On the other hand, it is  obvious that
$$
\pi_i(W^1+F^1(V^k(sl_3)))\subseteq W^1+F^1(V^k(sl_3)), \  i=2,3.
$$
So for $i=2,3$, we have
$$
\begin{array}{ll}
& \pi_i(\sum\limits_{i=2m}^{6m+1}a_ih_1(-1)^{i+1}h_2(-1)^{6m+1-i}(h_1+h_2)(-1){\bf 1})\\\\
=& c_i(\sum\limits_{i=2m}^{6m+1}a_ih_1(-1)^{i+1}h_2(-1)^{6m+1-i}(h_1+h_2)(-1){\bf 1}). 
\end{array}
$$
Notice that 
$$
\begin{array}{ll}
& \pi_2(\sum\limits_{i=2m}^{6m+1}a_ih_1(-1)^{i+1}h_2(-1)^{6m+1-i}(h_1+h_2)(-1){\bf 1})\\\\
=& \sum\limits_{i=2m}^{6m+1}a_i(h_1+h_2)(-1)^{i+1}(-h_2(-1))^{6m+1-i}h_1(-1){\bf 1}. 
\end{array}
$$
This means that $$\sum\limits_{i=2m}^{6m+1}a_ih_1(-1)^{i+1}h_2(-1)^{6m+1-i}(h_1+h_2)(-1){\bf 1}$$ has a factor $(h_1+h_2)(-1)^{2m+1}$.  Also
 $$
\begin{array}{ll}
& \pi_3(\sum\limits_{i=2m}^{6m+1}a_ih_1(-1)^{i+1}h_2(-1)^{6m+1-i}(h_1+h_2)(-1){\bf 1})\\\\
=& -\sum\limits_{i=2m}^{6m+1}a_ih_2(-1)^{i+1}h_1(-1)^{6m+1-i}(h_1+h_2)(-1){\bf 1}. 
\end{array}
$$
This implies that $\sum\limits_{i=2m}^{6m+1}a_ih_1(-1)^{i+1}h_2(-1)^{6m+1-i}(h_1+h_2)(-1){\bf 1}$ has a factor $h_2(-1)^{2m+1}$. Thus (\ref{e3.2.5c}) holds. 
\end{proof}
For the singular vector $v^2$, we can similarly obtain the following lemma.
\begin{lem}\label{l3.2.8c}  There exists a non-zero number $c$ such that 
	$$
	\begin{array}{ll}
	&cv^2-e_{\al_2}(-1)e_{\al_1+\al_2}(-1)h_1(-1)^{2m+1}h_2^{2m}(h_1+h_2)(-1)^{2m}{\bf 1}\\\\
	&-e_{\al_1}(-1)e_{\al_2}(-1)^2h_1(-1)^{2m}h_2(-1)^{2m}(h_1+h_2)(-1)^{2m}{\bf 1}\in W^1+F^1(V^k(sl_3)).
	\end{array}
	$$	
\end{lem}
We are now in a position to give the  second main result of this section. 
\begin{theorem}\label{t3.2.1} For $k=-3+\frac{2}{2m+1}$, $m\geq 1$, $X_{L_k(sl_3)}= \overline{G.({\mathbb C}^*(h_1-h_2)+f_{\theta})}$, where $\theta=\al_1+\al_2$.
\end{theorem}
\begin{proof} 
Notice that
$$\g=G. \C^*(h_1-h_2)\cup G\cdot (\C^*(h_1-h_2)+f_{\theta})\cup G\cdot ({\frak h}^{reg})\cup {\mathcal N},$$
and
$$
G. \C^*(h_1-h_2)\cup {\mathcal N}	\subseteq \overline{G\cdot (\C^*(h_1-h_2)+f_{\theta})}.
$$ 
By  Theorem \ref{tmain3.1},  the ideal generated by $v^1$ and $v^2$ is the maximal ideal of $V^k(sl_3)$. Recall that $W^1$ is the subspace of $V^k(sl_3)$ linearly spanned by monomials $z=z^{(+)}z^{(-)}z^{(0)}{\bf 1}$ such that $z\in F^0(V^k(sl_3))$ and $z^{(-)}\neq 1$. For $v\in V^k(sl_3)$, denote by $\bar{v}$ the image of $v$ in $R_{V^k(sl_3)}$. Then for any $w\in W^1$, the  
value of $\bar{w}$ at $h_1-h_2+f_{\theta}$ is zero. 
Let $h$ be a semisimple regular element in the Cartan subalgebra $\frak{h}=\C h_1+\C h_2$. Then 
all the $(h|h_1)$, $(h|h_2)$, and $(h|h_1+h_2)$ are non-zeros. Thus by Lemma \ref{l3.2.8}, 
\begin{equation}\label{e3.2.6c}
h\notin X_{L_k(sl_3)}.
\end{equation}
Then  by Lemmas \ref{l3.2.8}-\ref{l3.2.8c}, the zero locus of $\overline{{\mathcal I}_k}\subseteq R_{V^k(sl_3)}$ is exactly $\overline{G.(\C^*(h_1-h_2)+f_{\theta})}.$  It follows that
$$
X_{L_k(sl_3)}=\overline{G. (\C^*(h_1-h_2)+f_{\theta})}.
$$
\end{proof}
Recall that if $k\in\Z_{\geq 0}$, $X_{L_k(sl_3)}=\{0\}$ \cite{Zhu96,DM06, Ar12a},  $X_{L_k(sl_3)}=\overline{{\mathcal O}_{min}}$  if $k=-3+\frac{2m+1}{2}$, $m\geq 1$ \cite{Ar15a}, and $X_{L_k(sl_3)}={\mathcal N}$  if $k=-3$ or $k=-3+\frac{p}{q}$, $p,q\geq 3$, $(p,q)=1$ \cite{FF92}, \cite{Ar15a}.  Then we have 
\begin{cor} 
	\begin{enumerate}
		\item  For $k\in\C$, the associated variety $X_{L_k(sl_3)}$ of $L_k(sl_3)$ is one of the following:
		$$
		\g^*, \ \{0\}, \ \mathcal{N}, \ \overline{{\mathcal O}_{min}},  \ \overline{G. \C^*(h_1-h_2)}, \  \overline{G. (\C^*(h_1-h_2)+f_{\theta})}.
		$$		
		
		\item Let $h$ be a regular semi-simple vector of $sl_3$, then there is no $k\in\C$ such that 
		$X_{L_k(sl_3)}= \overline{G. \C^*h}$.
	\end{enumerate}
	\end{cor}
For $k=-3+\frac{2}{2m+1}$, $m\geq 0$,  let ${\mathcal I}_k$ be the maximal ideal of $V^k(sl_3)$, and ${\mathcal V}(I_k)$ be defined as in Section 2. We have the following conjecture.
\begin{conjecture}\label{conj1}
	$${\mathcal V}(I_k)=\{t\bar{\Lambda}_1-\frac{2i}{2m+1}\bar{\Lambda}_2, t\bar{\Lambda}_2-\frac{2i}{2m+1}\bar{\Lambda}_1,  t\bar{\Lambda}_1-(t+\frac{2i}{2m+1}+1)\bar{\Lambda}_2,  ~ t\in\C, i=0,1,\cdots, 2m\}.$$
\end{conjecture}
\begin{remark}
	If $m=0$, that is $k=-1$,  the conjecture is true by Proposition 5.5 of \cite{AP08}. When $m=1$, we could verify the conjecture  by computer programming.  We thank Libor K$\check{r}$i$\check{z}$ka  for pointing out the typos.
\end{remark}
\section{Varieties of simple affine $W$-algebras $W_k(sl_3,f)$}
\label{sec:W-algebras} 
We first have the following results from \cite{Ar15b} and \cite{Ara13}.
\begin{theorem}[\cite{Ar15b}] Let $\g$ be a simple complex Lie algebra, 
	 $f$  a regular  nilpotent element of $\g$,  and $k=-h^{\vee}+\frac{p}{q}$ a non-degenerate admissible number. Then the simple $W$-algebra $W_k(\g, f)$ is rational and lisse. 
\end{theorem}
\begin{theorem}[\cite{Ara13}]
		Let $f$ be a minimal nilpotent vector of $sl_3$ and $k=-3+\frac{2m+1}{2}$, $m\geq 1$. Then the simple $W$-algebra $W_k(sl_3)$ is rational and lisse. 
		\end{theorem}
Let $k=-1$ and $f$ a minimal nilpotent element of $sl_3$. It was proved in  \cite{AKMPP18} that  $H^0_{DS,f}(L_{-1}(sl_3))=W_{-1}(sl_3,f)$, and  is isomorphic to the rank one Heisenberg algebra $M(1)$. So  its associated variety  is one-dimensional. 

\begin{theorem}\label{thm5.1}
	 Let $k=-3+\frac{2}{2m+1}$, $m\geq 1$,  and $f$ a minimal nilpotent vector of $sl_3$, then 
	 \begin{enumerate}
	 	\item 
	 	\begin{equation}\label{w-5.3}
	 		X_{W_{k}(sl_3,f)}=\left\{\left[\begin{array}{ccc} a & b & 3(\mu^2-a^2)\\
	 			0 & -2a & c\\
	 			1 & 0 & a\end{array}\right]|~ \mu, a,b,c\in\C, \ bc=2(a-\mu)(2a+\mu)^2 \right\}.
	 	\end{equation}	
	 	In particular, 
	 	\begin{equation}\label{w-5.4}
	 		\dim X_{W_{k}(sl_3,f)}=3.
	 	\end{equation}
	 	
	 	\item $W_k(sl_3,f)$ is not quasi-lisse.
	 \end{enumerate}
	 
	\end{theorem}
\begin{proof} By Theorem \ref{Ar16},  for a minimal vector $f\in sl_3$, and $k+3=\frac{2}{2m+1}, m\in\Z_{\geq 0}$, 
	\begin{equation}\label{w-5.7c}
	H^0_{DS,f}(L_{k}(sl_3))=W_{k}(sl_3, f).  
	\end{equation}	
	By Theorem \ref{Th:W-algebra-variety} and (\ref{w-5.7c}),
	$$
	X_{W_k(sl_3,f)}=X_{L_{k}(sl_3)}\cap \mathscr{S}_{f},
	$$
for $k+3=\frac{2}{2m+1}$, $m\in\Z_{\geq 0}$. We may assume that $f=f_{\theta}$. Then
$$
\mathscr{S}_{f}=f_{\theta}+\g^{e_{\theta}},
$$	
and 
\begin{equation}\label{w-5.5}
\g^{e_{\theta}}=\C(\al_1-\al_2)\oplus \C e_{\theta}\oplus \C e_{\al_1}\oplus e_{\al_2}.
\end{equation}
Then we have
	\begin{equation}\label{w-5.6}
f_{\theta}+\g^{e_{\theta}}=\left\{\left[\begin{array}{ccc} a & b & d\\
0 & -2a & c\\
1 & 0 & a\end{array}\right]| ~a,b,c,d\in\C \right\}.
\end{equation}
It is obvious that for  $\mu\in\C$, 
$$
A_{\mu}=\left[\begin{array}{ccc} \mu & 0 & 0\\
0 & -2\mu& 0\\
1 & 0 & \mu\end{array}\right]\in  \overline{G\cdot (\C^*(h_1-h_2)+f_{\theta})}.
$$
In general,  non-semisimple
$$B=\left[\begin{array}{ccc} a & b & d\\
	0 & -2a & c\\
	1 & 0 & a\end{array}\right]\in \overline{G\cdot (\C^*(h_1-h_2)+f_{\theta})}\cap \mathscr{S}_{f_{\theta}}
	$$
 if and only if $B$  is similar to $A_{\mu}$ for some $\mu\in\C$. Then we deduce that
	$$
	d=3(\mu^2-a^2), 
	$$
	and 
	$$
	bc=2(a-\mu)(\mu+2a)^2.
	$$	
Thus (\ref{w-5.3}) and (\ref{w-5.4}) hold. 

\smallskip
By Theorem 8.2 of \cite{AMP23}, $L(t\bar{\Lambda}_1)$ and $L(t\bar{\Lambda}_2)$, $t\in\Z_{\geq 0}$, are irreducible ordinary modules of $L_k(sl_3)$. This  means that $L_k(sl_3)$ has infinitely many irreducible ordinary modules. 
 Then by Theorem 5.4 of \cite{AKR24}, 
for $t\in\Z_{\geq 0}$, the minimal quantum hamiltonian reductions ${\mathcal H}_{j_{t+1, 1}\Delta_{t+1,1}}$ and  ${\mathcal H}_{j_{1,t+1}\Delta_{1,t+1}}$
of $L(t\bar{\Lambda}_1)$ and $L(t\bar{\Lambda}_2)$ are  irreducible ordinary $W_k(sl_3,f)$-modules,  
where
$$
j_{t+1,1}=\frac{t}{3}, \  j_{1,t+1}=\frac{2t}{3}, \   \Delta_{t+1,1}=\frac{t^2+3t}{3(k+3)}-\frac{2t}{3}, \  \Delta_{1,t+1}=\frac{t^2+3t}{3(k+3)}-\frac{t}{3}.
$$
This means that $W_k(sl_3,f)$ has infinitely many irreducible modules. By Theorem \ref{ArK18}, $W_{k}(sl_3,f)$ is not quasi-lisse.
\end{proof}
We now assume that $f$ is a regular nilpotent vector of $sl_3$.  We may assume that 
$$
f=\left[\begin{array}{ccc} 0 & 0& 0\\
2 & 0 & 0\\
0 & 2 & 0\end{array}\right], \  e=\left[\begin{array}{ccc} 0 & 1 & 0\\
0 & 0 & 1\\
0 & 0 & 0\end{array}\right], \ h=\left[\begin{array}{ccc} 2 & 0 & 0\\
0 & 0 & 0\\
0 & 0 & -2\end{array}\right].
$$
Then
$$
\g^e=\C e\oplus \C e_{\theta},
$$
and 
	\begin{equation}\label{w-5.7}
f+\g^{e}=\left\{\left[\begin{array}{ccc} 0 & a & b\\
2 & 0 & a\\
0 & 2 & 0\end{array}\right]| ~a,b\in\C \right\}.
\end{equation}
Denote  the associated  universal affine $W$-algebras and the simple quotients by $W^k(sl_3)$ and $W_k(sl_3)$, respectively.  
\begin{theorem}\label{thm5.2} Let $k+3=\frac{2}{2m+1}$ or $\frac{2m+1}{2}$, then
	$$
	X_{W_k(sl_3)}=\left\{\left[\begin{array}{ccc} 0 & \frac{3}{4}\mu^2 & \frac{1}{2}\mu^3\\
	2 & 0 & \frac{3}{4}\mu^2\\
	0 & 2 & 0\end{array}\right]| ~\mu\in\C\right\}.
	$$
In particular, $\dim 	X_{W_k(sl_3)}=1$.  Furthermore, $W_k(sl_3)$ is not quasi-lisse.
\end{theorem}
\begin{proof} Since $k=-3+\frac{2}{2m+1}$, $m\geq 1$, it follows that
	$$
	(k\Lambda_0+\wr|\al)\notin \Z, \ {\rm for} \ \al\in\{-\al_1+\delta, -\al_2+\delta, -(\la_1+\al_2)+\delta, -(\al_1+\al_2)+2\delta\}.
	$$
	Then by Theorem 9.1.4 of \cite{Ar07}, 
	$$
	H^0_{DS,f}(L_{k}(sl_3))=W_{k}(sl_3).  
	$$
	Thus
	$$
	X_{W_k(sl_3))}=X_{L_{k}(sl_3)}\cap \mathscr{S}_{f}.
	$$
	By the Feigin-Frenkel Langlands duality \cite{FF92}, \cite{ACL19}, 
	$$
	W^k(sl_n)\cong W^{k'}(sl_n), 
	$$
	for $k+n=\frac{p}{q}$ and $k'+3=\frac{q}{p}$. In particular,  for $k+3=\frac{2}{2m+1}$ and $k'+3=\frac{2m+1}{2}$,
	$$
	W_k(sl_3)\cong W_{k'}(sl_3).
	$$
	Let 
	$$
	B=\left[\begin{array}{ccc} 0 & a& b\\
	2 & 0 & a\\
	0 & 2 & 0\end{array}\right]\in f+\g^e.
	$$
Then 
$$
f_{B}(\lambda)=|\lambda I-B|=\left|\begin{array}{ccc} \lambda & -a & -b\\
-2 & \lambda & -a\\
0 & -2 & \lambda\end{array}\right|=\lambda^3-4a\lambda-4b.
$$	
Notice that 
$$
X_{L_k(sl_3)}= \overline{G\cdot (\C^*(h_1-h_2)+f_{\theta})}.
$$	
So $B\in X_{W_k(sl_3)}$ if and only if $B$ is similar to
$$
A_{\mu}=\mu(h_1-h_2)+f_{\theta}=\left[\begin{array}{ccc} -\mu & 0 & 0\\
0 & 2\mu & 0\\
1 & 0 & -\mu\end{array}\right]
$$
for some $\mu\in\C$.  Then it is easy to deduce that $B\in X_{W_k(sl_3)}$ if and only if
 $$a=\frac{3}{4}\mu^2, \ b=\frac{1}{2}\mu^3.$$
 Thus
	$$
X_{W_k(sl_3)}=\left\{\left[\begin{array}{ccc} 0 & \frac{3}{4}\mu^2 & \frac{1}{2}\mu^3\\
2 & 0 & \frac{3}{4}\mu^2\\
0 & 2 & 0\end{array}\right]| ~\mu\in\C \right\}.
$$

Let $L(t\bar{\Lambda}_1)$ and $L(t\bar{\Lambda}_2)$, $t\in\Z_{\geq 0}$ be the irreducible ordinary modules of $L_k(sl_3)$ as in the proof of Theorem \ref{thm5.1}, then by  Theorem 5.4 of \cite{AKR24}, 
for $t\in\Z_{\geq 0}$, the principal  quantum hamiltonian reductions ${\mathcal W}_{h_{t+1, 1}w_{t+1,1}}$ and  ${\mathcal W}_{h_{1,t+1}w_{1,t+1}}$
of $L(t\bar{\Lambda}_1)$ and $L(t\bar{\Lambda}_2)$ are  irreducible ordinary $W_k(sl_3)$-modules,  
where
$$
h_{t+1,1}=\frac{t^2+3t}{3(k+3)}=h_{1,t+1},$$
and
$$ w_{t+1,1}=-\frac{2t\sqrt{3}}{3(k+3)^3}(\frac{2t}{3}-k-2)(\frac{t}{3}-k-2)=-w_{1,t+1}.
$$
It follows  from Theorem \ref{ArK18} that $W_k(sl_3)$ is not quasi-lisse. 
\end{proof}

\section{Declarations}
	
	\begin{enumerate}
		\item This work is supported by Natural Science Foundation of China (Grant number: 12171312);
		
		\item The authors have no relevant financial or non-financial interests to disclose;
		
		\item All datasets supporting the analysis and conclusions of the paper are publicly available at the time of publication.
		
		\end{enumerate}

\bibliographystyle{alpha}
\newcommand{\etalchar}[1]{$^{#1}$}

\end{document}